\documentstyle[12pt]{article}
\makeatletter
\@addtoreset{equation}{section}

\makeatother
\def\T{\displaystyle \mathop{\bf T}}
\def\R{\displaystyle \mathop{\bf R}}

\def\G1{\displaystyle \mathop{G}}
\def\e1{\displaystyle \mathop{e}}
\def\u1{\displaystyle \mathop{u}}

\def\T{\displaystyle \mathop{T}}
\def\n1{\displaystyle \mathop{\nu}}
\def\ps1{\displaystyle \mathop{\Psi}}

\def\p1{\displaystyle \mathop{p}}
\def\J1{\displaystyle \mathop{J}}

\def\O1{\displaystyle \mathop{O}}

\def\hG{\displaystyle \mathop{\hat{G}}}

\def\hgam{\displaystyle \mathop{\hat{\gamma}}}
\def\hpr{\displaystyle \mathop{\hat{\partial}}}
\def\Dlt{\displaystyle \mathop{\Delta}}
\def\hN{\displaystyle \mathop{\hat{N}}}
\def\he{\displaystyle \mathop{\hat{e}}}
\def\hp1{\displaystyle \mathop{\hat{p}}}
\def\hps{\displaystyle \mathop{\hat{\Psi}}}
\def\bhp{\displaystyle \mathop{\bar{\hat{\Psi}}}}
\def\bp{\displaystyle \mathop{\bar{\Psi}}}

\def\S{\displaystyle \sum}
\def\IIn{\displaystyle \int}
\def\FFr{\displaystyle \frac}
\def\Lm{\displaystyle \lim}
\def\pd{\displaystyle \prod}
\def\hb{\displaystyle \mathop{\hat{b}}}
\def\F{\displaystyle \mathop{F}}
\def\sig1{\displaystyle \mathop{\sigma}}

\def\h1{\displaystyle \mathop{h}}
\def\H1{\displaystyle \mathop{H}}
\def\A1{\displaystyle \mathop{A}}

\def\D1{\displaystyle \mathop{D}}
\def\D1{\displaystyle \mathop{D}}
\def\g1{\displaystyle \mathop{g}}

\begin{document}
\begin{center}
{\Large {\bf Mathematical Background of Formalism }}\\
{\Large {\bf of Operator Manifold}}
\vskip 0.5truecm
{\normalsize G.T.Ter-Kazarian\footnote
{E-mail address:gago@bao.sci.am}}\\
{\small Byurakan Astrophysical Observatory, Armenia 378433}\\
{\small September 3, 1997}\\
\end{center}
\begin{abstract}
The analysis of mathematical structure of the method of 
operator manifold [1] guides our discussion. The nature of operator manifold
provides its elements with both quantum field and 
differential geometry aspects, a detailed study of which is a subject of  
present paper. It yields a quantization of geometry differing in 
principle from all earlier suggested schemes. This formalism was made complete 
by construction of state wave functions and calculation of matrix elements
of geometric objects. While, it has been shown that the matrix element 
of any geometric object of operator manifold gave rise to corresponding 
geometric object of wave manifold.
\end{abstract}
\vskip 0.5truecm
\section {Introduction}
\label {int}
In spite of considerable progress achieved over the entire subsequent 
period in the study of fundamental constituents of matter and forces,
the physical theory is still far from being complete and could not be 
regarded as the final word in particle physics, since many fundamental 
questions have yet to be answered. The absence of the vital theory 
which will be able to solve the crucial problems of particle physics 
imperatively stimulates a search for general constructive principles.
Effecting a reconciliation in [1] we are led to consider the theory
exploring the query of origin of geometry, some basic concepts of particle
physics and also four major principles of Relativity, Quantum, Gauge 
and Color Confinement. We have finally arrived at the scenario
a whole idea of which comes to following: the geometry,
quarks with various quantum numbers, internal symmetries and etc.
also those four principles, as it was proven, are derivative
and come into being simultaneously.\\
In [1] we elaborated a new mathematical framework, which is a still 
wider generalization of the method of secondary quantization
with appropriate expansion over the geometric objects.
The formalism of {\em operator manifold} ensued, which is a 
guiding formalism framing our approach yielding a quantization of
geometry differing in principle from all earlier suggested schemes. The
nature of operator manifold provides its elements with both field and 
geometric aspects. To save writing, in [1] they have been discussed briefly,
leaving a more profound mathematical study for separate treatment. 
As far as suggested theory involves a drastic revision of our ideas 
of geometry, basic concepts and principles of particle physics, in order 
to be more consistent and convinced in correctness of drawn statements 
it will be advantageous to fill this shortage and verify them by further 
exploration of structure of the theory. 
There is an attempt to supply the mathematical background and to trace
some of the major currents of thoughts under its view-point.
In view of all this,
the detailed study of different aspects surely an important 
subject for present article. It can be regarded as a mathematical addition 
to [1]. The article is organized as follows: To begin with our task will
be to make few preliminary
remarks on the structure of manifold $G(2.2.3)$, on the base of which
we define general conceptual part of formalism of operator manifold
$\hat{G}(2.2.3)$. In the aftermath, we analyze the mathematical structure
of the method of $\hat{G}(2.2.3)$ in both aspects of quantum field theory
and differential geometry.
This is not a final report on a closed subject, but it is hoped that suggested 
theory will serve as useful introduction and that it will thereby add a
knowledge of method of operator manifold in quantization of geometry.

\vskip 0.5truecm

\section {Operator-Vector and Co-Vector Fields}
\label{Vec}

Before embarking on the main strategy of discussion, as a starting point, 
just a very brief recapitulation of major properties of structure of manifold
$G(2.2.3)$ [1-6]. 
Let the maximal curve $\lambda(t):R^{1}\rightarrow G(2.2.3)$ passed through
the point $p=\lambda(0)\in G(2.2.3)$ with tangent vector 
$\left.{\bf A}\right|_{\lambda(t)}$, where the $G(2.2.3)$ is 12-dimensional
smooth differentiable manifold. The set $\{\zeta^{(\lambda,\mu,\alpha)} \}$
$(\lambda,\mu=12; \quad \alpha=1,2,3)$ are local coordinates in open 
neighborhood of $p\in \cal U$, namely the curve $\lambda(t)$ has the
coordinates $\zeta^{(\lambda,\mu,\alpha)}(t)$. For any other point
$q\in G(2.2.3)$ it can be found open neighborhood $\cal U$ and 
$\varepsilon > 0$ such that ${\bf A}$ gives the diffeomorphism
$f_{t}:{\cal U}\rightarrow G(2.2.3)$ at $|t| <\varepsilon$. 
The 12-dimensional smooth vector field ${\bf A}_{p} ={\bf A}({\bf \zeta})$ 
belongs to the section of tangent bundle ${\bf T}_{p}$ at the point
$p({\bf \zeta})$. So define the differential $d\,A^{t}_{p}$ of the flux
$A^{t}_{p}:G(2.2.3)\rightarrow G(2.2.3)$ at the point $p({\bf \zeta})\in
G(2.2.3)$ with the field of velocities ${\bf A}({\bf \zeta})$, then
one gets $d\,A^{t}_{p}:{\bf T}_{p}\rightarrow R$, namely the one-parameter
group of diffeomorphismes $A^{t}$ given for the maximal curve 
${\bf \zeta}(t)$ passing through point $p$ and ${\bf \zeta}(0)=
{\bf \zeta}_{p}, \quad \dot{{\bf \zeta}}(0)={\bf A}_{p}$
\begin{equation}
\label{R21}
d\,A^{t}_{p}({\bf A})=\left.\FFr{d}{d\,t}\right|_{t=0}A^{t}\left(
{\bf \zeta}(t)\right) = {\bf A}_{p}({\bf \zeta}).
\end{equation}
Hence, 
\begin{center}
$d\,A^{t}:{\bf T}(G(2.2.3))\rightarrow R$ $\left( 
{\bf T}(G(2.2.3)) = \displaystyle \mathop{\bigcup}_{p({\bf \zeta})}
{\bf T}_{p}\right)$.
\end{center}
Let the 
$\{e_{(\lambda,\mu,\alpha)}=O_{\lambda,\mu}\otimes
\sigma_{\alpha}\} \subset G(2.2.3)$
is some set of linear independent $12$ unit vectors at the point $p$,
provided with the linear unit bi-pseudo-vectors $\{ O_{\lambda,\mu} \}$:
$<O_{\lambda,\mu},O_{\tau,\nu}>={}^{*}\delta_{\lambda,\tau}
{}^{*}\delta_{\mu,\nu}$
serving as the basis for tangent vectors of 
$2 \times 2$ dimensional linear bi-pseudo-space 
${}^{*}G(2.2)$, ${}^{*}\delta_{\lambda,\tau}=1-\delta_{\lambda,\tau}$, 
where $\delta$ is Kronecker symbol, ; 
and the ordinary unit vectors $ \sigma_{\alpha}$
implying $<\sigma_{\alpha}, \sigma_{\beta}>= \delta_{\alpha\beta}$.
Henceforth we always let the first two subscripts in the parentheses
specify the bi-pseudo-vector components, while the third refers to the 
ordinary-vector components.
The metric on $G(2.2.3)$ is bilinear, local, symmetric and positive 
defined reflection of vector fields of sections ${\bf T}$ of tangent 
bundle of $G(2.2.3)$, namely  
$\hat{\bf g}:{\bf T}_{p}\otimes {\bf T}_{p}\rightarrow C^{\infty}(G(2.2.3))$ 
a section of conjugate vector bundle $S^{2}{\bf T}$ (symmetric part of 
tensor degree) with components in basis $\{e_{(\lambda,\mu,\alpha)}\}$
\begin{equation}
\label{R22}
g_{(\lambda,\mu,\alpha)(\tau,\nu,\beta)}=g(e_{(\lambda,\mu,\alpha)},
e_{(\tau,\nu,\beta)})=g(e_{(\tau,\nu,\beta)},e_{(\lambda,\mu,\alpha)}).
\end{equation}
Any vector ${\bf A}_{p}\in{\bf T}_{p}$ reads ${\bf A}=
e_{(\lambda,\mu,\alpha)}A^{(\lambda,\mu,\alpha)}$, provided with
components $A^{(\lambda,\mu,\alpha)}$ in the basis
$\{e_{(\lambda,\mu,\alpha)}\}$.
Except where stated otherwise, here as usual,
the double occurrence of the dummy indices will be taken
to denote a summation extended over their all values.
In holonomic coordinate basis $\left(\partial/\partial\,
\zeta^{(\lambda,\mu,\alpha)}\right)_{p}$ one gets 
$A^{(\lambda,\mu,\alpha)}=\left.\FFr{d\,
\zeta^{(\lambda,\mu,\alpha)}}{d\,t}\right|_{p}$ and
$\hat{g}=g_{(\lambda,\mu,\alpha)(\tau,\nu,\beta)}d\zeta^
{(\lambda,\mu,\alpha)}\otimes d\zeta^{(\tau,\nu,\beta)}$.
The manifold $G(2.2.3)$ decomposes as follows:
$G(2.2.3)={}^{*}{\bf R}^{2}
\otimes {}^{*}{\bf R}^{2} \otimes {\bf R}^{3}=\G1_{\eta}(2.3)\oplus
\G1_{u}(2.3)=\displaystyle\sum^{2}_{\lambda,\mu=1} 
\oplus {\bf R}^{3}_{\lambda \mu}=
{\R_{x}}^{3}\oplus 
{\T_{x}}^{3}\oplus
{\R_{u}}^{3}\oplus 
{\T_{u}}^{3}$ with corresponding 
basis vectors  
${\e1_{i}}^{0}_{(\lambda\alpha)}={\O1_{i}}_{\lambda}\otimes 
\sigma_{\alpha}
\subset \G1_{i}(2.3)$ $(i=\eta, u)$ of tangent sections, where 
${\O1_{i}}_{+}=
\displaystyle \frac{1}{\sqrt{2}}(O_{1,1} +\varepsilon_{i} O_{2,1})$,
${\O1_{i}}_{-}=
\displaystyle \frac{1}{\sqrt{2}}(O_{1,2} +\varepsilon_{i} O_{2,2})$,
$\varepsilon_{\eta}=1$, $\varepsilon_{u}=-1$. There up on
$<{\O1_{i}}_{\lambda},{\O1_{i}}_{\tau}>=
\varepsilon_{i}\delta_{ij}{}^{*}\delta_{\lambda \tau}$.
The positive metric forms are defined
in manifolds $\G1_{i} (2.3)$   
$\eta^{2}=\eta_{(\lambda\alpha)}\eta^{(\lambda\alpha)}
\in \G1_{\eta}(2.3), \quad
u^{2}=u_{(\lambda\alpha)} u^{(\lambda\alpha)}\in \G1_{u}(2.3),$
where
\begin{center}
$
\begin{array}{l}
\eta^{(+ \alpha)}= \FFr{1}{\sqrt{2}}(\zeta^{(1,1,\alpha)}+
\zeta^{(2,1,\alpha)}),\quad  
\eta^{(- \alpha)}= \FFr{1}{\sqrt{2}}(\zeta^{(1,2,\alpha)}+
\zeta^{(2,2,\alpha)}), \\ 
u^{(+ \alpha)}= \FFr{1}{\sqrt{2}}
(\zeta^{(1,1,\alpha)}-\zeta^{(2,1,\alpha)}),\quad
u^{(- \alpha)}= \FFr{1}{\sqrt{2}}
(\zeta^{(1,2,\alpha)}-\zeta^{(2,2,\alpha)}).
\end{array}
$
\end{center}
The $G(2.3)$ decomposes into three-dimensional
ordinary $({\bf R}^{3})$ and time $({\bf T}^{3})$ flat spaces $G(2.3)=
{\bf R}^{3}\oplus {\bf T}^{3}$ with signatures $sgn({\bf R}^{3})=(+++)$ 
and $sgn({\bf T}^{3})=(---)$. Since all directions in ${\bf T}^{3}$ are
equivalent, then by notion {\em time} one implies the projection of
time-coordinate on fixed arbitrary universal direction in ${\bf T}^{3}$.
By the reduction ${\bf T}^{3}\rightarrow  T^{1}$ the transition
$G(2.3)\rightarrow M^{4}={\bf R}^{3}\oplus T^{1}$ may be performed
whenever it will be needed. 
At this point we cut short a discussion of structure of 
$G(2.2.3)$ and refer to [3,4] for details.\\
Next we proceed to preliminary definitions of the elements of operator
manifold, which will be further discussed and made complete in due
course of exposition of the formalism describing the processes of 
creation and annihilation of geometric objects. This formalism is 
analogous to the method of secondary quantization describing the 
processes of creation and annihilation of particles in the 
configuration space of occupation numbers, but  
it will be an appropriate expansion over the geometric objects.
Adjusting to fit a conventional notations, below we change 
the order of vector's and co-vector's indices used in [1] into the opposite.
So, the state of $\zeta$-type quantum is describe by means of function
$\Phi^{(\lambda,\mu,\alpha)}(\zeta)$ belonging to the ring of functions
of $C^{\infty}$-class: 
\begin{center}
$
\begin{array}{l}
\Phi^{(1,1,\alpha)}=\FFr{1}{\sqrt{2}}\left({\ps1_{\eta} }
^{(+\alpha)}+{\ps1_{u}}^{(+\alpha)}\right),\quad 
\Phi^{(1,2,\alpha)}=\FFr{1}{\sqrt{2}}\left({\ps1_{\eta} }
^{(-\alpha)}+{\ps1_{u}}^{(-\alpha)}\right), \\ 
\Phi^{(2,1,\alpha)}=\FFr{1}{\sqrt{2}}\left({\ps1_{\eta} }
^{(+\alpha)}+{\ps1_{u}}^{(+\alpha)}\right), \quad
\Phi^{(2,2,\alpha)}=\FFr{1}{\sqrt{2}}\left({\ps1_{\eta} }
^{(-\alpha)}+{\ps1_{u}}^{(-\alpha)}\right),
\end{array}
$
\end{center}
provided by the functions of $\eta$- and $u$-type quanta defined on
the manifolds $\G1_{i}(2.3)$
\begin{equation}
\label {R23}
\begin{array}{l}
{\ps1_{\eta} }^{(\pm\alpha)}(\eta,p_{\eta})=\eta^{(\pm\alpha)}
{\ps1_{\eta} }^{\pm}(\eta,p_{\eta}),\\ 
{\ps1_{u}}^{(\pm\alpha)}(u,p_{u})=u^{(\pm\alpha)}
{\ps1_{u}}^{\pm}(u,p_{u}).
\end{array}
\end{equation}
Here it was assumed that the probability of finding the quantum in the 
state with fixed coordinate ($\eta$ or $u$) and momentum 
($p_{\eta}$ or $p_{u}$) is determined by the square of its state 
wave function ${\ps1_{\eta} }_{\pm}(\eta,p_{\eta}),$ or 
${\ps1_{u}}_{\pm}(u,p_{u})$. This provides a simple intuitive meaning
of state functions, where
the $6$-vectors of coordinates - $\eta,u$, and
momenta  - $p_{\eta},p_{u}$ respectively are
$\eta= {\e1_{\eta}}_{(\lambda\alpha)}^{0}\eta^{(\lambda\alpha)}, \quad
p_{\eta}= {\e1_{\eta}}_{(\lambda\alpha)}^{0}
{\p1_{\eta}}^{(\lambda\alpha)}, \quad
u= {\e1_{u}}_{(\lambda\alpha)}^{0}u^{(\lambda\alpha)}, \quad
p_{u}= {\e1_{u}}_{(\lambda\alpha)}^{0}
{\p1_{u}}^{(\lambda\alpha)}$. \\
Being confronted by the problem of quantization of geometry, we
first deal with a substitution of the basis elements by the 
corresponding operators of creation and annihilation of quanta acting 
in the configuration space of occupation numbers.
Instead of pseudo-vectors $O_{\lambda}$ we introduce the following 
operators supplied by additional index ($r$) referring to the quantum 
numbers of corresponding state
\begin{equation}
\label {R24}
\begin{array}{ll}
\hat{O}^{r}_{1}=O^{r}_{1}\alpha_{1},\quad
\hat{O}^{r}_{2}=O^{r}_{2}{\alpha}_{2},\quad
\hat{O}_{r}^{\lambda}={}^{*}\delta^{\lambda\mu}\hat{O}^{r}_{\mu}=
{(\hat{O}^{r}_{\lambda})}^{+},\\
\{ \hat{O}^{r}_{\lambda},\hat{O}^{r'}_{\tau} \}=
\delta_{rr'}{}^{*}\delta_{\lambda\tau}I_{2}, 
\quad
<{O}^{r}_{\lambda},{O}^{r'}_{\tau}>= \delta_{rr'}{}^{*}\delta_{\lambda\tau},
\qquad I_{2}=\left( \begin{array}{cc}
1 \quad 0 \\
0 \quad 1
\end{array}
\right).
\end{array}
\end{equation}
The matrices ${\alpha}_{\lambda}$ satisfy the conditions
\begin{equation}
\label {R25}
\begin{array}{lr}
{\alpha}^{\lambda}={}^{*}\delta^{\lambda\mu}
{\alpha}_{\mu}={({\alpha}_{\lambda})}^{+},\qquad
\{ {\alpha}_{\lambda},{\alpha}_{\tau} \}={}^{*}\delta_
{\lambda\tau}I_{2}.
\end{array}
\end{equation}
For example, they may be in the form
${\alpha}_{1}=\left( \begin{array}{cc}
0 \quad 1 \\
0 \quad 0
\end{array}
\right), \quad
{\alpha}_{2}=\left( \begin{array}{cc}
0 \quad 0 \\
1 \quad 0
\end{array}
\right).$
This forms the starting point for quantization.
Creation operator $\hat{O}^{r}_{1}$ generates one-occupied state
$\mid 1>_{(0)}\equiv\mid 0,\ldots,1,\ldots>_{(0)}$ and the basis vector
$O^{r}_{1}$ with the quantum number $r$
right through acting on non-occupied vacuum state
$\mid 0>_{(0)}\equiv  \mid 0,0,\ldots>_{(0)}$:
\begin{equation}
\label{R25}
\hat{O}^{r}_{1}\mid 0>_{(0)}={O}^{r}_{1}\mid 1>_{(0)}.
\end{equation}
Accordingly, the action of annihilation operator $\hat{O}^{r}_{2}$
on one-occupied state yields the vacuum state and the basis vector
$O^{r}_{2}$
\begin{equation}
\label{R26}
\hat{O}^{r}_{2}\mid 1>_{(0)}={O}^{r}_{2}\mid 0>_{(0)}.
\end{equation}
So define 
$\hat{O}^{r}_{1}\mid 1>_{(0)}=0,\quad \hat{O}^{r}_{2}\mid 0>_{(0)}=$0.
The matrix realization of the states $\mid 0>_{0}$ and $\mid 1>_{0}$,
for instance, may be as follows:
$\mid 0>_{0}\equiv\chi_{1}=\left( \begin{array}{c}
0 \\
1
\end{array}
\right), \quad
\mid 1>_{0}\equiv\chi_{2}=\left( \begin{array}{c}
1  \\
0
\end{array}
\right).$
The operator of occupation number is
$\hat{N}^{(0)}_{r}=\hat{O}^{r}_{1}\hat{O}^{r}_{2}$,
with the expectation values implying Pauli's exclusion principle
${}_{(0)}<0\mid\hat{N}^{(0)}_{r}\mid 0>_{(0)}=0, \qquad 
{}_{(0)}<1\mid\hat{N}^{(0)}_{r}\mid 1>_{(0)}=1$.
The vacuum state reads
$\chi_{0}\equiv\mid 0>_{(0)}=\displaystyle\prod_{r=1}^{N}(\chi_{1})_{r}$.
With this final detail cared for one-occupied state takes the form
$\chi_{r'}\equiv\mid 1>_{(0)}=(\chi_{2})_{r'}\displaystyle\prod_{r\neq r'}
(\chi_{1})_{r}$.
Continuing along this line, instead of ordinary vectors we introduce
the operators
$\hat{\sigma}^{r}_{\alpha}\equiv\delta_{\alpha\beta\gamma}
\sigma^{r}_{\beta}\widetilde{\sigma}_{\gamma}$,
where $\widetilde{\sigma}_{\gamma}$ are Pauli's matrices, and
\begin{equation}
\label{R27}
<\sigma_{\alpha}^{r},\sigma_{\beta}^{r'}>=\delta_{rr'}\delta_{\alpha\beta},
\quad
\hat{\sigma}^{\alpha}_{r}=\delta^{\alpha\beta}\hat{\sigma}^{r}_{\beta}=
{(\hat{\sigma}_{\alpha}^{r})}^{+}=\hat{\sigma}_{\alpha}^{r},
\quad
\{\hat{\sigma}_{\alpha}^{r},\hat{\sigma}_{\beta}^{r'}\}=2
\delta_{rr'}\delta_{\alpha\beta}I_{2}.
\end{equation}
For the vacuum state $\mid 0>_{(\sigma)}\equiv{\varphi}_{1(\alpha)}$
and one-occupied state $\mid 1_{(\alpha)}>_{(\sigma)}
\equiv{\varphi}_{2(\alpha)}$
we make use of matrix realization
${\varphi}_{1(\alpha)}\equiv\chi_{1}, \quad
{\varphi}_{2(1)}=\left( \begin{array}{c}
1 \\
0
\end{array}
\right), \quad
{\varphi}_{2(2)}=\left( \begin{array}{c}
-i \\
\hspace{0.25cm} 0
\end{array}
\right), \quad
{\varphi}_{2(3)}=\left( \begin{array}{c}
\hspace{0.25cm} 0 \\
-1
\end{array}
\right).$
Then
\begin{equation}
\label{R28}
{\hat{\sigma}}_{\alpha}^{r}\varphi_{1(\alpha)}=\sigma_
{\alpha}^{r}\varphi_{2(\alpha)}=(\sigma_{\alpha}^{r}\widetilde{\sigma}_
{\alpha})\varphi_{1(\alpha)},
\quad
{\hat{\sigma}}_{\alpha}^{r}\varphi_{2(\alpha)}=\sigma_
{\alpha}^{r}\varphi_{1(\alpha)}=(\sigma_{\alpha}^{r}\widetilde{\sigma}_
{\alpha})\varphi_{2(\alpha)}.
\end{equation}
Hence, the single eigen-value
$(\sigma_{\alpha}^{r}\widetilde{\sigma}_{\alpha})$
has associated with it quite different $\varphi_{\lambda(\alpha)}$.
The eigen-value is degenerated with 
degeneracy degree equal 2. Due to it, along many quantum numbers $r$
there is also the quantum number of the spin $\vec{\sigma}$ with the values
$\sigma_{3}=\FFr{1}{2}s\quad (s=\pm1)$.
This rule for spin quantum number is not without an important reason. 
The argument for this conclusion
is compulsory suggested by the properties of operators 
$\hat{\sigma}^{r}_{\alpha}$.
As it was seen in [1], this consequently {\em gives rise to spin of
particle}.\\
One-occupied state reads
$\varphi_{r'(\alpha)}={(\varphi_{2(\alpha)})}_{r'}\displaystyle
\prod_{r\neq r'}{(\chi_{1})}_{r}$.
Next we introduce the operators
\begin{equation}
\label{R29}
{\hat{\gamma}}^{r}_{(\lambda,\mu,\alpha)}\equiv{\hat{O}}^{r_{1}}_{\lambda}
\otimes{\hat{O}}^{r_{2}}_{\mu}\otimes{\hat{\sigma}}^{r_{3}}_{\alpha},
\quad
{\hat{\gamma}}_{r}^{(\lambda,\mu,\alpha)}\equiv{\hat{O}}_{r_{1}}^{\lambda}
\otimes{\hat{O}}_{r_{2}}^{\mu}\otimes{\hat{\sigma}}_{r_{3}}^{\alpha}=
{}^{*}\delta^{\lambda\tau}{}^{*}\delta^{\mu\nu}\delta^{\alpha\beta} 
{\hat{\gamma}}^{r}_{(\tau,\nu,\beta)},
\end{equation}
and also the state vector
\begin{equation}
\label{R210}
\chi_{\lambda,\mu,\tau(\alpha)}\equiv\mid\lambda,\mu,\tau(\alpha)>=
\chi_{\lambda}\otimes\chi_{\mu}\otimes\varphi_{\tau(\alpha)}, 
\end{equation}
where $\lambda,\mu,\tau,\nu=
1,2;\quad \alpha,\beta=1,2,3$ and $r\equiv (r_{1},r_{2},r_{3})$.
Hence
${\hat{\gamma}}^{r}_{(\lambda,\mu,\alpha)}\chi_{\tau,\nu,\delta(\beta)}=
({\hat{O}}^{r_{1}}_{\lambda}\chi_{\tau})\otimes
({\hat{O}}^{r_{2}}_{\mu}\chi_{\nu})\otimes({\hat{\sigma}}^{r_{3}}_{\alpha}
\varphi_{\delta(\beta)})$.
Omitting the two-valuedness of state vector we apply
$\mid\lambda,\tau,\delta(\beta)>\equiv\mid\lambda,\tau>$,
and the same time remember that always the summation 
must be extended over the double degeneracy of the spin states $(s=\pm 1)$.\\
With this final detail cared for one infers the explicit forms of corresponding
matrix elements:
\begin{equation}
\label{R211}
\begin{array}{l}
<2,2\mid{\hat{\gamma}}^{r}_{(1,1,\alpha)}\mid 1,1>=e^{r}_{(1,1,\alpha)},\quad
<1,1\mid{\hat{\gamma}}_{r}^{(1,1,\alpha)}\mid 2,2>=e_{r}^{(1,1,\alpha)},\\
<2,1\mid{\hat{\gamma}}^{r}_{(1,2,\alpha)}\mid 1,2>=e^{r}_{(1,2,\alpha)},\quad
<1,2\mid{\hat{\gamma}}_{r}^{(1,2,\alpha)}\mid 2,1>=e_{r}^{(1,2,\alpha)},\\
<1,2\mid{\hat{\gamma}}^{r}_{(2,1,\alpha)}\mid 2,1>=e^{r}_{(2,1,\alpha)},\quad
<2,1\mid{\hat{\gamma}}_{r}^{(2,1,\alpha)}\mid 1,2>=e_{r}^{(2,1,\alpha)},\\
<1,1\mid{\hat{\gamma}}^{r}_{(2,2,\alpha)}\mid 2,2>=e^{r}_{(2,2,\alpha)},\quad
<2,2\mid{\hat{\gamma}}_{r}^{(2,2,\alpha)}\mid 1,1>=e_{r}^{(2,2,\alpha)}.
\end{array}
\end{equation}
The operators of occupation numbers are
\begin{equation}
\label{R212}
\begin{array}{l}
{\hN_{1}}^{rr'}_{\alpha\beta}=
{\hat{\gamma}}^{r}_{(1,1,\alpha)}{\hat{\gamma}}^{r'}_{(2,2,\beta)}=
{\hN_{1}}_{rr'}^{\alpha\beta}=
{\hat{\gamma}}_{r}^{(2,2,\alpha)}{\hat{\gamma}}_{r'}^{(1,1,\beta)}, \\
{\hN_{2}}^{rr'}_{\alpha\beta}=
{\hat{\gamma}}^{r}_{(2,1,\alpha)}{\hat{\gamma}}^{r'}_{(1,2,\beta)}=
{\hN_{2}}_{rr'}^{\alpha\beta}=
{\hat{\gamma}}_{r}^{(1,2,\alpha)}{\hat{\gamma}}_{r'}^{(2,1,\beta)},
\end{array}
\end{equation}
with the expectation values implying Pauli's exclusion principle
\begin{equation}
\label{R214}
\begin{array}{ll}
<2,2\mid{\hN_{1}}_{rr'}^{\alpha\beta}\mid 2,2>=\delta_{rr'}
\delta_{\alpha\beta},
\quad
<1,2\mid{\hN_{2}}_{rr'}^{\alpha\beta}\mid 1,2>=\delta_{rr'}\delta_
{\alpha\beta},\\
<1,1\mid{\hN_{1}}_{rr'}^{\alpha\beta}\mid 1,1>=0,\qquad
<2,1\mid{\hN_{2}}_{rr'}^{\alpha\beta}\mid 2,1>=0.
\end{array}
\end{equation}
The set of operators $\{{\hat{\gamma}}^{r}_{(\lambda,\mu,\alpha)}\}$
is the basis for operator-vectors 
$\hat{\Phi}(\zeta)={\hat{\gamma}}^{r}_{(\lambda,\mu,\alpha)}
\Phi_{r}^{(\lambda,\mu,\alpha)}(\zeta)$,
but a set of operators $\{{\hat{\gamma}}_{r}^{(\lambda,\mu,\alpha)}\}$
is a dual basis for
operator-co-vectors
$\bar{\hat{\Phi}}(\zeta)={\hat{\gamma}}_{r}^{(\lambda,\mu,\alpha)}
\Phi^{r}_{(\lambda,\mu,\alpha)}(\zeta)$, where
$\Phi^{r}_{(\lambda,\mu,\alpha)}(\zeta)=
{\bar{\Phi}}_{r}^{(\lambda,\mu,\alpha)}(\zeta)$ (charge-conjugated).
One easily gets
\begin{equation}
\label{R215}
\begin{array}{l}
<2,2\mid\hat{\Phi}(\zeta)\bar{\hat{\Phi}}(\zeta)\mid 2,2>=
\Phi_{r}^{(1,1,\alpha)}(\zeta)\Phi^{r}_{(1,1,\alpha)}(\zeta),\\
<2,1\mid\hat{\Phi}(\zeta)\bar{\hat{\Phi}}(\zeta)\mid 2,1>=
\Phi_{r}^{(1,2,\alpha)}(\zeta)\Phi^{r}_{(1,2,\alpha)}(\zeta),\\
<1,2\mid\hat{\Phi}(\zeta)\bar{\hat{\Phi}}(\zeta)\mid 1,2>=
\Phi_{r}^{(2,1,\alpha)}(\zeta)\Phi^{r}_{(2,1,\alpha)}(\zeta),\\
<1,1\mid\hat{\Phi}(\zeta)\bar{\hat{\Phi}}(\zeta)\mid 1,1>=
\Phi_{r}^{(2,2,\alpha)}(\zeta)\Phi^{r}_{(2,2,\alpha)}(\zeta).
\end{array}
\end{equation}
Introducing the  state vectors
\begin{equation}
\label{R216}
\begin{array}{l}
\chi^{0}(\nu_{1},\nu_{2},\nu_{3},\nu_{4})=
\mid 1,1>^{\nu_{1}}\cdot\mid 1,2>^{\nu_{2}}\cdot
\mid 2,1>^{\nu_{3}}\cdot\mid 2,2>^{\nu_{4}},\\
\\
\nu_{i}= \left\{ \begin{array}{ll}
                   1   & \mbox{if $\nu=\nu_{i}$}\quad  \mbox{for some $i$,} \\
                   0   & \mbox{otherwise},
                   \end{array}
\right. \\
\\
\mid\chi_{-}(\lambda)>=\left\{ \begin{array}{ll}
                   \chi^{0}(1,0,0,0)   & \lambda=1, \\
                   \chi^{0}(0,0,1,0)   & \lambda=2,
                   \end{array}
\right. \quad
\mid\chi_{+}(\lambda)>=\left\{ \begin{array}{ll}
                   \chi^{0}(0,0,0,1)   & \lambda=1, \\
                   \chi^{0}(0,1,0,0)   & \lambda=2,
                   \end{array}
\right.
\end{array}
\end{equation}
provided
$$
\begin{array}{ll}
<\lambda,\mu,\mid\tau,\nu>=\delta_{\lambda\tau}\delta_{\mu\nu}, \quad
<\chi_{\pm}\mid A\mid \chi_{\mp}>=
\S_{\lambda}<\chi_{\pm}(\lambda)\mid A\mid \chi_{\mp}(\lambda)>,\\
<\chi_{\pm}(\lambda)\mid\chi_{\pm}(\mu)>=\delta_{\lambda\mu}, \quad
<\chi_{\pm}\mid A\mid \chi_{\pm}>=
\S_{\lambda}<\chi_{\pm}(\lambda)\mid A\mid \chi_{\pm}(\lambda)>,\\
<\chi_{\pm}(\lambda)\mid\chi_{\mp}(\mu)>=0,
\end{array}
$$
we get the matrix elements as follows:
\begin{equation}
\label{R217}
\begin{array}{l}
<\chi_{+}\mid\hat{\Phi}(\zeta)\bar{\hat{\Phi}}(\zeta)\mid\chi_{+}>
\equiv\Phi^{2}_{+}(\zeta)=\quad
\Phi_{r}^{(1,1,\alpha)}(\zeta)\Phi^{r}_{(1,1,\alpha)}(\zeta)+
\Phi_{r}^{(2,1,\alpha)}(\zeta)\Phi^{r}_{(2,1,\alpha)}(\zeta),\\
<\chi_{-}\mid\hat{\Phi}(\zeta)\bar{\hat{\Phi}}(\zeta)\mid\chi_{-}>
\equiv\Phi^{2}_{-}(\zeta)=
\Phi_{r}^{(2,2,\alpha)}(\zeta)\Phi^{r}_{(2,2,\alpha)}(\zeta)+
\Phi_{r}^{(1,2,\alpha)}(\zeta)\Phi^{r}_{(1,2,\alpha)}(\zeta).
\end{array}
\end{equation}
The basis $\{{\hat{\gamma}}^{r}_{(\lambda,\mu,\alpha)}\}$
decomposes into
$\{ {\hgam_{i}}^{r}_{(\lambda\alpha)} \}\quad
(\lambda=\pm;\quad\alpha=1,2,3;\quad i=\eta,u).$
The latter reads in component form
${\hgam_{i}}^{r}_{(+\alpha)}=\FFr{1}{\sqrt{2}}
({\hgam}^{r}_{(1,1\alpha)}+\varepsilon_{i}
{\hgam}^{r}_{(2,1\alpha)}),
\quad
{\hgam_{i}}^{r}_{(-\alpha)}=\FFr{1}{\sqrt{2}}
({\hgam}^{r}_{(1,2\alpha)}+\varepsilon_{i}
{\hgam}^{r}_{(2,2\alpha)}).$ 
The expansions of operator-vectors $\hps_{i}\in\hG_{i}(2.3)$ and 
operator-co-vectors  $\bar{\hps_{i}}\in\hG_{i}(2.3)$ are written 
$\hps_{i}={\hgam_{i}}^{r}_{(\lambda\alpha)}{\ps1_{i}}_{r}^{(\lambda\alpha)},
\quad\bar{\hps_{i}}=
{\hgam_{i}}_{r}^{(\lambda\alpha)}{\ps1_{i}}^{r}_{(\lambda\alpha)},$
where the components ${\ps1_{\eta}}_{r}^{(\lambda\alpha)}(\eta)$ and
${\ps1_{u}}_{r}^{(\lambda\alpha)}(u)$ are in the form eq.(2.3), and
${\ps1_{i}}^{r}_{(\lambda\alpha)}=
{\bar{\ps1_{i}}}_{r}^{(\lambda\alpha)}$.
The operator of occupation number of $i$-type quantum takes the form
${\hN_{i}}_{rr'}^{\alpha\beta}=\varepsilon_{i}{\hgam_{i}}_{r}^{(-\alpha)}
{\hgam_{i}}_{r'}^{(+\beta)}$,
with corresponding expectation values
\begin{equation}
\label {R218}
<\chi_{-}\mid{\hN_{i}}_{rr'}^{\alpha\beta}\mid\chi_{-}>=0, 
\quad
<\chi_{+}\mid{\hN_{i}}_{rr'}^{\alpha\beta}\mid\chi_{+}>
=\varepsilon_{i}<{\e1_{i}}_{r}^{(-\alpha)},{\e1_{i}}_{r'}^{(+\beta)}>=
\delta_{rr'}\delta_{\alpha\beta}.
\end{equation}
Taking into account two-valuedness of degenerate spin states
it follows that the quanta are the fermions with
half-integral spins. That is, the functions ${\ps1_{\eta}}^
{\lambda}(\eta)$ and ${\ps1_{u}}^{\lambda}(u)$ may be regarded as the
Fermi fields of $\eta$- and $u$-type quanta.
Explicitly the matrix elements read
\begin{equation}
\label {R219}
\begin{array}{l}
\Phi^{2}_{+}(\zeta)=<\chi_{+}\mid\hps_{\eta}(\eta)\bhp_{\eta}(\eta)
+ \hps_{u}(u)\bhp_{u}(u)\mid\chi_{+}>=
\varepsilon_{i}{\ps1_{i}}^{(+\alpha)}{\ps1_{i}}_{(+\alpha)}=\\
={\ps1_{\eta}}^{(+\alpha)}(\eta){\bp_{\eta}}^{(+\alpha)}(\eta)
-{\ps1_{u}}^{(+\alpha)}(u)
{\bp_{u}}^{(+\alpha)}(u),\\
\Phi^{2}_{-}(\zeta)=<\chi_{-}\mid\hps_{\eta}(\eta)\bhp_{\eta}(\eta)
+ \hps_{u}(u)\bhp_{u}(u)\mid\chi_{-}>=
\varepsilon_{i}{\ps1_{i}}^{(-\alpha)}{\ps1_{i}}_{(-\alpha)}=\\
={\ps1_{\eta}}^{(-\alpha)}(\eta)
{\bp_{\eta}}^{(-\alpha)}(\eta)
-{\ps1_{u}}^{(-\alpha)}(u)
{\bp_{u}}^{(-\alpha)}(u).  
\end{array}
\end{equation}

\section {Operator Manifold $\hat{G}(2.2.3)$} 
\label {quant}
The field aspect will be the subject for discussion in this section. 
Our notation will be that of the textbook
by [7].
As far as a quantum may be regarded as a fermion field, its description
is provided by the theory, which is in 
close analogy to Dirac's conventional wave-mechanical theory of fermions
with spin $\vec{\FFr{1}{2}}$ treated in terms of manifold 
$G(2.2.3)$ [1]. The final formulation of quantum theory 
is equivalent to configuration space wave mechanics with antisymmetric state
functions.
To facilitate our approach it seems useful to present some
formal matters which one will have to know in order to understand the
structure of theory. Below we proceed with preliminary discussion.
The quantum field may be considered as
bi-spinor field
$\Psi(\zeta)$ defined on manifold $G(2.2.3)=\G1_{\eta}(2.3)
\oplus\G1_{u}(2.3)$:  $\Psi(\zeta)=\ps1_{\eta}(\eta)\ps1_{u}(u)$,
where the $\ps1_{i}$ is a bi-spinor defined on 
the manifold $\G1_{i}(2.3).$
The state of free quantum of $i$-type with definite values of link-momentum
$\p1_{i}$ and spin projection $s$ is described by means of plane waves,
respectively (in units $\hbar=1,\quad c=1$):
\begin{equation}   
\label{R31}
\begin{array}{l}
{\ps1_{\eta}}_{p_{\eta}}(\eta)={\left({\FFr{m}{E_{\eta}}}
\right)} ^{1/2}
\u1_{\eta}(p_{\eta},s)e^{-ip_{\eta}\eta}, 
\quad
{\ps1_{u}}_{p_{u}}(u)={\left( {\FFr{m}{E_{u}}}
\right)} ^{1/2}
\u1_{u}(p_{u},s)e^{-ip_{u}u},
\end{array}
\end{equation}
where
$E_{i}\equiv{\p1_{i}}_{0}={({\p1_{i}}_{0\alpha},{\p1_{i}}_{0\alpha})}^{1/2},
\quad
{\p1_{i}}_{0\alpha}=\FFr{1}{\sqrt{2}}({\p1_{i}}_{(+\alpha)}+
{\p1_{i}}_{(-\alpha)})$.
It is necessary to consider also the solutions of negative
frequencies
\begin{equation}   
\label{R32}
{\ps1_{\eta}}_{-p_{\eta}}(\eta)={\left({\FFr{m}{E_{\eta}}}
\right) }^{1/2}
\n1_{\eta}(p_{\eta},s)e^{ip_{\eta}\eta}, 
\quad
{\ps1_{u}}_{-p_{u}}(u)={\left( {\FFr{m}{E_{u}}}
\right) }^{1/2}
\n1_{u}(p_{u},s)e^{ip_{u}u},
\end{equation}
where
${\p1_{i}}_{\alpha}=\FFr{1}{\sqrt{2}}({\p1_{i}}_{(+\alpha)}-
{\p1_{i}}_{(-\alpha)}),
\quad
p^{2}_{\eta}=E^{2}_{\eta}-{\vec{p}}^{2}_{\eta}=p^{2}_{u}=
E^{2}_{u}-{\vec{p}}^{2}_{u}=m^{2}$.
For the spinors the useful relations of orthogonality and
completeness hold.
We make use of localized wave packets constructed by
means of superposition of plane wave solutions furnished by creation
and annihilation operators in agreement with Pauli's principle
\begin{equation}
\label{R33}
\hps_{\eta}(\eta)= \S_{\pm s}\int\frac{d^{3}p_{\eta}}{{(2\pi)}^{3/2}}
\hps_{\eta}(p_{\eta},s,\eta),
\quad
\hps_{u}(u)= \S_{\pm s}\int\frac{d^{3}p_{u}}{{(2\pi)}^{3/2}}
\hps_{u}(p_{u},s,u),
\end{equation}
where, as usual, it is denoted
\begin{equation}   
\label{R34}
\begin{array}{ll}
\hps_{\eta}(p_{\eta},s,\eta)={\hgam_{\eta}}_{(\lambda\alpha)}(p_{\eta},s)
{\ps1_{\eta}}^{(\lambda\alpha)}(p_{\eta},s,\eta),
\quad
\bhp_{\eta}(p_{\eta},s,\eta)={\hgam_{\eta}}^{(\lambda\alpha)}(p_{\eta},s)
{\ps1_{\eta}}_{(\lambda\alpha)}(p_{\eta},s,\eta),\\
\hps_{u}(p_{u},s,u)={\hgam_{u}}_{(\lambda\alpha)}(p_{u},s)
{\ps1_{u}}^{(\lambda\alpha)}(p_{u},s,u),
\quad
\bhp_{u}(p_{u},s,u)={\hgam_{u}}^{(\lambda\alpha)}(p_{u},s)
{\ps1_{u}}_{(\lambda\alpha)}(p_{u},s,u).
\end{array}
\end{equation}
One has
${\ps1_{\eta}}^{(\pm\alpha)}=\eta^{(\pm\alpha)}{\ps1_{\eta}}^{\pm},\quad
{\ps1_{u}}^{(\pm\alpha)}=u^{(\pm\alpha)}{\ps1_{u}}^{\pm}, 
\quad
{\ps1_{\eta}}_{(\pm\alpha)}=\eta_{(\pm\alpha)}{\ps1_{\eta}}_{\pm},\quad
{\ps1_{u}}_{(\pm\alpha)}=u_{(\pm\alpha)}{\ps1_{u}}_{\pm}, $
provided
${\ps1_{i}}^{+}\equiv{\ps1_{i}}_{p_{i}},\quad{\ps1_{i}}^{-}\equiv
{\ps1_{i}}_{-p_{i}},\quad{\ps1_{i}}_{\lambda}={\bar{\ps1_{i}}}^{\lambda},
\quad{\ps1_{i}}_{(\lambda\alpha)}={\bar{\ps1_{i}}}^{(\lambda\alpha)}$.
A closer examination of the properties of the matrix elements of 
the anticommutators of expansion coefficients shows that
\begin{equation}   
\label{R35}
\begin{array}{l}
<\chi_{-}\mid \{ {\hgam_{i}}^{(+\alpha)}(p_{i},s),
{\hgam_{j}}_{(+\beta)}(p'_{j},s')\}\mid\chi_{-}>=\\
=<{\he_{i}}^{(+\alpha)}(p_{i},s),{\he_{j}}_{(+\beta)}(p'_{j},s')>=
\varepsilon_{i}\delta_{ij}\delta_{ss'}\delta_{\alpha\beta}\delta^{(3)}
({\vec{p}}_{i}- {\vec{p'}}_{i}),\\
<\chi_{+}\mid \{ {\hgam_{j}}^{(-\beta)}(p'_{j},s'),
{\hgam_{i}}_{(-\alpha)}(p_{i},s)\}\mid\chi_{+}>=\\
=<{\he_{j}}^{(-\beta)}(p'_{j},s'),{\he_{i}}_{(-\alpha)}(p_{i},s)>=
\varepsilon_{i}\delta_{ij}\delta_{ss'}\delta_{\alpha\beta}\delta^{(3)}
({\vec{p}}_{i}- {\vec{p'}}_{i}).
\end{array}
\end{equation}
We may also consider analogical wave packets of operator-vector
fields $\hat{\Phi}(\zeta)$ and $\bar{\hat{\Phi}}(\zeta)$.
While, explicitly the matrix element of anticommutator reads
\begin{equation}
\label{R36}
\begin{array}{l}
<\chi_{\pm}\mid\{ {\hat{\gamma}}^{(\lambda,\mu,\alpha)}(p,s),
{\hat{\gamma}}_{(\tau,\nu,\beta)}(p',s')\}\mid\chi_{\pm}>=\\
=<e^{(\lambda,\mu,\alpha)}(p,s),e_{(\tau,\nu,\beta)}(p',s')>=
\delta_{ss'}\delta^{\lambda}_{\tau}\delta^{\nu}_{\mu}\delta^{\alpha}_{\beta}
\delta^{(3)}(\vec{p}-\vec{p'}).
\end{array}
\end{equation}
Thus
\begin{equation}
\label{R37}
\begin{array}{l}
\S_{\lambda=\pm}<\chi_{\lambda}\mid\hat{\Phi}(\zeta)
\bar{\hat{\Phi}}(\zeta)\mid
\chi_{\lambda}>= 
\S_{\lambda=\pm}<\chi_{\lambda}\mid
\bar{\hat{\Phi}}(\zeta)\hat{\Phi}(\zeta)\mid\chi_{\lambda}>= \\
-i\Lm_{\zeta\rightarrow\zeta'}(\zeta\zeta')\G1_{\zeta}(\zeta-\zeta')=
-i\left[ \Lm_{\eta\rightarrow\eta'}(\eta\eta')\G1_{\eta}(\eta-\eta')-
\Lm_{u\rightarrow u'}(uu')\G1_{u}(u-u')\right],
\end{array}
\end{equation}
where the Green's functions are used
\begin{equation}
\label{R38}
\G1_{\eta}(\eta-{\eta}')=-(i\hpr_{\eta}+m)\Dlt_{\eta}(\eta-{\eta}'),
\quad
\G1_{u}(u-u')=-(i\hpr_{u}+m)\Dlt_{u}(u-u'),
\end{equation}
provided with the invariant singular functions 
$\Dlt_{\eta}(\eta-{\eta}')$ and $\Dlt_{u}(u-u')$. Meanwhile a second
trend emerged as the quantum theory situation corresponding
to simultaneously presence of many identical quanta. 
To describe the $n$-quanta fermion system by means of quantum field theory, 
it will be advantageous to make use of convenient method of constructing 
the state vector of physical system by proceeding from the vacuum state
as a very point of origin.
Exploiting the whole advantage of it, a particular emphasis will be
placed just on the fact that,
certainly, if there are $n$ identical quanta with coordinates
$\zeta_{1},\zeta_{2},...,\zeta_{n}$ the antisymmetrical state function
$\Phi$ will be a function of all of them and presents the system of
$n$ fermions $\Phi(\zeta_{1},\zeta_{2},...,\zeta_{n})$, which implies
the Fermi-Dirac statistics.\\
It was assumed that the $i$-th quantum is found in the state $r_{i}$ 
with the field function $\Phi_{r_{i}}^{(\lambda_{i},\mu_{i},\alpha_{i})}=
\zeta_{r_{i}}^{(\lambda_{i},\mu_{i},\alpha_{i})}
\Phi_{r_{i}}^{\lambda_{i},\mu_{i}}(\zeta_{r_{i}})$ and made use of 
following notation: 
$\zeta_{r_{i}^{\lambda_{i},\mu_{i}}}=$\\
$\S^{3}_{\alpha_{i}=1}
e^{r_{i}}_{(\lambda_{i},\mu_{i},\alpha_{i})}
\zeta_{r_{i}}^{(\lambda_{i},\mu_{i},\alpha_{i})}$;
$\zeta_{r_{i}}=\S^{2}_{\lambda_{i},\mu_{i}=1}
\zeta_{r_{i}^{\lambda_{i},\mu_{i}}}\in \widetilde{\cal U}_{r_{i}}$, 
the $\widetilde{\cal U}
_{r_{i}}$
is the open neighborhood of the point $\zeta_{r_{i}};$ the $r_{i}$ 
implies a set $\left(r_{i}^{11},r_{i}^{12},r_{i}^{21},r_{i}^{22}\right)$.
Let the ${\cal H}^{(1)}$ is a Hilbert space used for quantum mechanical 
description of one particle, namely ${\cal H}^{(1)}$ is a finite or infinite
dimensional complex space, provided with scalar product $(\Phi,\Psi)$
being linear with respect to $\Psi$ and antilinear to $\Phi$.  
The ${\cal H}^{(1)}$ is complete in sense of norm 
$|\Phi|=(\Phi,\Phi)^{1/2}$,
i.e. each fundamental sequence $\{\Phi_{n}\}$ of vectors of ${\cal H}^{(1)}$
converged by norm in ${\cal H}^{(1)}$. 
One-particle state function is written 
$\Phi^{(1)}_{r_{i}}=\displaystyle \prod^{2}_{\lambda_{i},\mu_{i}=1}
\Phi^{(1)}_{r_{i}^{\lambda_{i},\mu_{i}}}\in {\cal H}^{(1)}_{r_{i}}$, where
${\cal H}^{(1)}_{r_{i}}=
\displaystyle \prod^{2}_{\lambda_{i},\mu_{i}=1}\otimes
{\cal H}^{(1)}_{r_{i}^{\lambda_{i},\mu_{i}}}$. So define
\begin{equation}
\label{R39}
\widetilde{\Phi}^{(1)}=\zeta_{i}\Phi^{(1)}_{r_{i}}\in 
\widetilde{G}^{(1)}_{r_{i}}=\widetilde{\cal U}^{(1)}_{r_{i}}\otimes
{\cal H}^{(1)}_{r_{i}}.
\end{equation}
For description of n-particle system we introduce Hilbert space
\begin{equation}
\label{R319}
{\bar{\cal H}}^{(n)}_{(r_{1},\ldots,r_{n})}=
{\cal H}^{(1)}_{r_{1}}\otimes\cdots\otimes
{\cal H}^{(1)}_{r_{n}}
\end{equation}
by considering all sequences
\begin{equation}
\label{R310}
\Phi^{(n)}_{(r_{1},\ldots,r_{n})}=
\{\Phi^{(1)}_{r_{1}},\ldots,\Phi^{(1)}_{r_{n}}\}=
\Phi^{(1)}_{r_{1}}\otimes \cdots \otimes \Phi^{(1)}_{r_{n}},
\end{equation}
where $\Phi^{(1)}_{r_{i}}\in {\cal H}^{(1)}_{r_{i}}$, provided, 
as usual, with the scalar product
\begin{equation}
\label{R312}
(\Phi^{(n)}_{(r_{1},\ldots,r_{n})},\Psi^{(n)}_{(r_{1},\ldots,r_{n})})=
\prod^{n}_{i=1}(\Phi^{(1)}_{r_{i}},\Psi^{(1)}_{r_{i}}).
\end{equation}
Obtaining the space eq.(3.9) we consider the space 
${\cal H}^{(n)}_{(r_{1},\ldots,r_{n})}$ of all limited linear combinations 
of eq.(3.10) and continue by linearity the scalar product eq.(3.12)
on ${\cal H}^{(n)}_{(r_{1},ldots,r_{n})}$. The wave function
$\Phi^{(n)}_{(r_{1},\ldots,r_{n})} 
\in {\cal H}^{(n)}_{(r_{1},\ldots,r_{n})}$ 
must be antisymmetrized over its arguments. So, we ought to distinguish the
antisymmetric part ${}^{A}{\bar{\cal H}}^{(n)}$ of Hilbert space
${\bar{\cal H}}^{(n)}$ by considering the functions
\begin{equation}
\label{R313}
{}^{A}\Phi^{(n)}_{(r_{1},\ldots,r_{n})}=\FFr{1}{\sqrt{n!}}\S_{\sigma\in S(n)}
sgn(\sigma)\Phi^{(n)}_{\sigma(r_{1},\ldots,r_{n})}.
\end{equation}
The summation is extended over all permutations of indices
$(r_{1}^{\lambda\mu},\ldots, r_{n}^{\lambda\mu})$ of the integers
$1,2,\ldots,n,$ whereas the antisymmetrical eigen-functions are sums
of the same terms with alternating signs in dependence of a parity
$sgn(\sigma)$ of transposition. In the aftermath, one continues 
by linearity on ${\cal H}^{(n)}$ the reflection 
$\Phi^{(n)}\rightarrow {}^{A}\Phi^{(n)}$ , which is limited and 
allowed the expansion by linearity on ${}^{A}{\bar{\cal H}}^{(n)}$.
According to conventional definition, the region of values of this
reflection is a ${}^{A}{\bar{\cal H}}^{(n)}$, namely the antisymmetrized
tensor product of $n$ identical samples of ${}{\cal H}^{(1)}$.
We introduce
\begin{equation}
\label{R314}
\begin{array}{l}
{}^{A}{\widetilde{\Phi}}^{(n)}_{(r_{1},\ldots,r_{n})}=
\FFr{1}{\sqrt{n!}}\S_{\sigma\in S(n)}
sgn(\sigma){\widetilde {\Phi}}^{(n)}_{\sigma(r_{1},\ldots,r_{n})}=\\
=\FFr{1}{\sqrt{n!}}\S_{\sigma\in S(n)}
sgn(\sigma){\widetilde{\Phi}}^{(1)}_{r_{1}}\otimes \cdots \otimes 
{\widetilde{\Phi}}^{(1)}_{r_{n}} \in 
{}^{A}{\widetilde{G}}^{(n)}_{(r_{1},\ldots,r_{n})}=
\widetilde{\cal U}^{(n)}_{(r_{1},\ldots,r_{n})}\otimes
{}^{A}{\hat{\cal H}}^{(n)}_{(r_{1},\ldots,r_{n})}.
\end{array}
\end{equation}
and constructing 12-dimensional wave manifold $\widetilde{G}(2.2.3)$ 
consider a set ${}^{A}\widetilde{\cal F}$ of all sequences
\begin{equation}
\label{R315}
{}^{A}{\widetilde{\Phi}}=\{ {}^{A}{\widetilde{\Phi}}^{(0)},
{}^{A}{\widetilde{\Phi}}^{(1)}\ldots ,
{}^{A}{\widetilde{\Phi}}^{(n)}\ldots \},
\end{equation}
with a finite number of non-zero elements. Therewith, the set 
${}^{A}\cal F$
\begin{equation}
\label{R316}
{}^{A}\Phi=\{ {}^{A}\Phi^{(0)},
{}^{A}\Phi^{(1)}\ldots ,
{}^{A}\Phi^{(n)}\ldots \}
\end{equation}
is provided with the structure of Hilbert sub-space by introducing
the composition rules
\begin{equation}
\label{R317}
\begin{array}{l}
{}^{A}(\lambda\Phi +\mu\Psi)^{(n)}=
\lambda {}^{A}\Phi^{(n)} + \mu {}^{A}\Psi^{(n)}, \quad 
\forall \lambda, \mu \in C,\\
({}^{A}\Phi,{}^{A}\Psi)=\S^{\infty}_{n=0}({}^{A}\Phi^{(n)},{}^{A}\Psi^{(n)}).
\end{array}
\end{equation}
In the sequel, the wave manifold $\widetilde{G}(2.2.3)$ stems from the
${}^{A}\widetilde{\cal F}$ by making expansion by metric 
induced as a scalar product in ${}^{A}\cal F$
\begin{equation}
\label{R318}
\widetilde{G}(2.2.3)=\S^{\infty}_{n=0}\widetilde{G}^{(n)}=
\S^{\infty}_{n=0}\left(\widetilde{\cal U}^{(n)}\otimes {}^{A}
\bar{\cal H}^{(n)}\right).
\end{equation}
In general, the creation ${\hat{\gamma}}_{r}^{(\lambda,\mu,\alpha)}$ and
annihilation ${\hat{\gamma}}^{r}_{(\lambda,\mu,\alpha)}$ operators for each
${}{\cal H}^{(1)}$ can be defined as follows.
One ought to modify the operators eq.(2.10) in order to provide
an anticommutation being valid in both cases acting on the same as well as
different states:
\begin{equation}
\label{R319}
{\hat{\gamma}}_{(\lambda,\mu,\alpha)}^{r}\Rightarrow
{\hat{\gamma}}_{(\lambda,\mu,\alpha)}^{r}\eta_{r}^{\lambda\mu},
\quad
{\hat{\gamma}}^{(\lambda,\mu,\alpha)}_{r}\Rightarrow
\eta_{r}^{\lambda\mu}{\hat{\gamma}}^{(\lambda,\mu,\alpha)}_{r}=
{\hat{\gamma}}^{(\lambda,\mu,\alpha)}_{r}\eta_{r}^{\lambda\mu},
\quad
{(\eta_{r}^{\lambda\mu})}^{+}=\eta_{r}^{\lambda\mu},
\end{equation}
for fixed $\lambda,\mu,\alpha$,
where $\eta_{r}^{\lambda\mu}$ is a diagonal operator in the space 
of occupation numbers.
Therewith, at $r_{i}<r_{j}$ one gets
\begin{equation}
\label{R320}
{\hat{\gamma}}_{(\lambda,\mu,\alpha)}^{r_{i}}\eta_{r_{j}}^{\lambda\mu}=
-\eta_{r_{j}}^{\lambda\mu}{\hat{\gamma}}_{(\lambda,\mu,\alpha)}^{r_{i}},
\quad
{\hat{\gamma}}_{(\lambda,\mu,\alpha)}^{r_{j}}\eta_{r_{i}}^{\lambda\mu}=
\eta_{r_{i}}^{\lambda\mu}{\hat{\gamma}}_{(\lambda,\mu,\alpha)}^{r_{j}}.
\end{equation}
The operators of corresponding occupation numbers, for fixed 
$\lambda,\mu,\alpha$, are
\begin{equation}
\label{R321}
{\hat{N}}_{r}^{\lambda\mu}={\hat{\gamma}}_{(\lambda,\mu,\alpha)}^{r}
{\hat{\gamma}}^{(\lambda,\mu,\alpha)}_{r}={\hat{N}}_{r}^{\lambda}\otimes
{\hat{N}}_{r}^{\mu},
\end{equation}
where we make use of (see eq.(2.5))
\begin{equation}
\label{R322}
{\hat{N}}^{\lambda}_{r}=\left\{ \begin{array}{l}
\hat{O}^{r}_{1}\hat{O}_{r}^{1}=(\alpha_{1}\alpha_{2})_{r}= \left( 
\begin{array}{ll}
1 &0 \\
0 &0
\end{array} \right)_{r}, \\
\\
\hat{O}^{r}_{2}\hat{O}_{r}^{2}=(\alpha_{2}\alpha_{1})_{r}= \left( 
\begin{array}{ll}
0 &0 \\
0 &1
\end{array} \right)_{r}.
\end{array} \right.
\end{equation}
As far as diagonal operators $(1-2{\hat{N}}_{r}^{\lambda\mu})$ 
anticommute with the
${\hat{\gamma}}_{(\lambda,\mu,\alpha)}^{r}$, then
\begin{equation}
\label{R323}
\eta_{r_{i}}^{\lambda\mu}=\prod_{r =1}^{r_{i}-1}(1-2{\hat{N}}_{r}
^{\lambda\mu}).
\end{equation}
Combining eq.(3.20)-eq.(3.23), the explicit forms read
\begin{equation}
\label{R324}
\begin{array}{ll}
\eta_{r_{i}}^{11}=\pd_{r =1}^{r_{i}-1}\left( \begin{array}{rrrr}
-1 &0 &0 &0\\
 0 &1 &0 &0\\
 0 &0 &1 &0\\
 0 &0 &0 &1
\end{array}
\right)_{r}, 
\quad
\eta_{r_{i}}^{21}=\pd_{r =1}^{r_{i}-1}\left( \begin{array}{rrrr}
 1 &0 &0  &0\\
 0 &1 &0  &0\\
 0 &0 &-1 &0\\
 0 &0 &0  &1
\end{array}
\right)_{r}, \\
\\
\eta_{r_{i}}^{12}=\pd_{r =1}^{r_{i}-1}\left( \begin{array}{rrrr}
 1 &0  &0 &0\\
 0 &-1 &0 &0\\
 0 &0  &1 &0\\
 0 &0  &0 &1
\end{array}
\right)_{r}, 
\quad
\eta_{r_{i}}^{22}=\pd_{r =1}^{r_{i}-1}\left( \begin{array}{rrrr}
 1 &0 &0 &0\\
 0 &1 &0 &0\\
 0 &0 &1 &0\\
 0 &0 &0 &-1
\end{array}
\right)_{r}, 
\end{array}
\end{equation}
provided
\begin{equation}
\label{R325}
\begin{array}{l}
\eta_{r_{i}}^{11}\Phi(n_{1},\ldots,n_{N};0;0;0)=\pd_{r =1}^{r_{i}-1}
(-1)^{n_{r}}\Phi(n_{1},\ldots,n_{N};0;0;0),\\
\\
\eta_{r_{i}}^{12}\Phi(0;m_{1},\ldots,n_{M};0;0)=\pd_{r =1}^{r_{i}-1}
(-1)^{m_{r}}\Phi(0;m_{1},\ldots,m_{M};0;0),\\
\\
\eta_{r_{i}}^{21}\Phi(0;0;q_{1},\ldots,q_{Q};0)=\pd_{r =1}^{r_{i}-1}
(-1)^{q_{r}}\Phi(0;0;q_{1},\ldots,q_{Q};0),\\
\\
\eta_{r_{i}}^{22}\Phi(0;0;0;t_{1},\ldots,t_{T})=\pd_{r =1}^{r_{i}-1}
(-1)^{t_{r}}\Phi(0;0;0;t_{1},\ldots,t_{T}),
\end{array}
\end{equation}
Here the occupation numbers $n_{r}(m_{r},q_{r},t_{r})$ are introduced, 
which refer to the $r$-th states corresponding to operators
${\hat{\gamma}}_{(1,1,\alpha)}^{r}({\hat{\gamma}}_{(1,2,\alpha)}^{r},
{\hat{\gamma}}_{(2,1,\alpha)}^{r},{\hat{\gamma}}_{(2,2,\alpha)}^{r})$
either empty ($n_{r},\ldots,t_{r}=0$) or occupied
($n_{r},\ldots,t_{r}=1$).
To save writing we abbreviate the modified operators by the
same symbols. For example,
the creation operator ${\hat{\gamma}}^{(\lambda,\mu,\alpha)}_{r_{i}}$ by
acting on free state $\mid 0>_{r_{i}}$ yields the one-occupied state
$\mid 1>_{r_{i}}$ with the phase $+$ or $-$ depending of parity of the number
of quanta in the states $ r < r_{i}$.
Modified operators satisfy the same anticommutation relations
of the operators eq.(2.10). 
It is convenient to make use of notation
${\hat{\gamma}}^{(\lambda,\mu,\alpha)}_{r}\equiv
{e}^{(\lambda,\mu,\alpha)}_{r}{\hat{b}}^{\lambda\mu}_{(r\alpha)},
\quad
{\hat{\gamma}}_{(\lambda,\mu,\alpha)}^{r}\equiv
{e}_{(\lambda,\mu,\alpha)}^{r}{\hat{b}}_{\lambda\mu}^{(r\alpha)}$,
and abbreviate the pair of indices $(r\alpha)$ by the single symbol $r$. 
Then, for each $\Phi \in {}^{A}{\cal H}^{(n)}$
\begin{equation}
\label{R326}
\Phi=\FFr{1}{\sqrt{n!}}\S_{\sigma\in S(n)}
sgn(\sigma)\Phi^{(n)}_{\sigma(r_{1},\ldots,r_{n})},
\end{equation}
where $\Phi^{(n)}_{\sigma(r_{1},\ldots,r_{n})}=
\Phi^{(1)}_{\sigma{(1)}}\otimes\cdots\otimes \Phi^{(1)}_{\sigma{(n)}}$,
and for any vector $f \in {\cal H}^{(1)}$, the operators
$\hat{b}(f)$ $(\hat{b}_{\lambda\mu}^{r}(f))$ and
$\hat{b}^{*}(f)$ $(\hat{b}^{\lambda\mu}_{r}(f))$ imply
\begin{equation}
\label{R327}
\begin{array}{l}
\hat{b}(f)\Phi=\FFr{1}{\sqrt{(n-1)!}}\S_{\sigma\in S(n)}
sgn(\sigma)\left( f\,\Phi^{(1)}_{\sigma{(1)}} \right)
\Phi^{(1)}_{\sigma{(2)}}\otimes\cdots\otimes \Phi^{(1)}_{\sigma{(n)}},\\
\hat{b}^{*}(f)\Phi=\FFr{1}{\sqrt{(n+1)!}}\S_{\sigma\in S(n+1)}
sgn(\sigma)\Phi^{(1)}_{\sigma{(0)}}\otimes
\Phi^{(1)}_{\sigma{(1)}}\otimes\cdots\otimes \Phi^{(1)}_{\sigma{(n)}},
\end{array}
\end{equation}
where $\Phi^{(1)}_{(0)}\equiv f$. In the aftermath, one continues the 
$\hat{b}(f)$ and $\hat{b}^{*}(f)$ by linearity to linear reflections, 
which will be denoted by the same symbols, acting from ${}^{A}{\cal H}^{(n)}$
into ${}^{A}{\cal H}^{(n-1)}$ or ${}^{A}{\cal H}^{(n+1)}$, respectively.
They were limited by the values $\sqrt{n}|f|$ and $\sqrt{(n+1)}|f|$.
So, they may be expanded by continuation up to the reflections acting 
from ${}^{A}{\bar{\cal H}}^{(n)}$
into ${}^{A}{\bar{\cal H}}^{(n-1)}$ or ${}^{A}{\bar{\cal H}}^{(n+1)}$.
Finally, they ought to be continued by linearity up to the linear
operators acting from ${}^{A}{\cal F}$ into ${}^{A}{\cal F}$.
They are defined on the same close region in
${}^{A}{\bar{\cal H}}^{(n)}$, namely in ${}^{A}{\cal F}$, which is 
invariant with respect to reflections
$\hat{b}(f)$ and $\hat{b}^{*}(f)$. Hence, at $f_{i}, g_{i}\in 
{\cal H}^{(1)}$ $(i=1,\ldots,n; j=1,\ldots,m)$ all polynomials over 
$\{\hat{b}^{*}(f_{i}) \}$ and $\{\hat{b}(g_{j}) \}$ are completely defined 
on ${}^{A}{\cal F}$. While
\begin{equation}
\label{R328}
\begin{array}{l}
<1,1\mid\{ {\hat{b}}^{11}_{r}(f),{\hat{b}}_{11}^{r'}(g)
\}\mid 1,1>=
<1,2\mid\{ {\hat{b}}^{12}_{r}(f),{\hat{b}}_{12}^{r'}(g)
\}\mid 1,2>=\\
= <2,1\mid\{ {\hat{b}}^{21}_{r}(f),{\hat{b}}_{(21}^{r'}(g)
\}\mid 2,1>=
<2,2\mid\{ {\hat{b}}^{22}_{r}(f),{\hat{b}}_{22}^{r'}(g)
\}\mid 2,2> = \delta^{r'}_{r}.
\end{array}
\end{equation}
The mean values $<\varphi; {\hat{b}}_{r}^{\lambda\mu}(f)
{\hat{b}}^{r}_{\lambda\mu}(f)>$ calculated at given $\lambda,\mu$ for
any element $\Phi \in {}^{A}{\cal F}$ equal to mean values of the symmetric 
operator of occupation number in terms of
${\hat{N}}_{\lambda\mu}^{r}={\hat{b}}_{r}^{\lambda\mu}(f)
{\hat{b}}^{r}_{\lambda\mu}(f)$, with a wave function $f$ in the 
state describing by $\Phi$. 
Here, as usual, it was denoted $<\varphi; A\Phi> = Tr \,
P_{\varphi}A = (\Phi,A\Phi)$ for each vector $\Phi \in {\cal H}$ with
$|\Phi |=1$, while the $P_{\varphi}$ is projecting operator onto one-
dimensional space $\{\lambda \Phi\left.\right| \lambda\in C\}$ generated by
$\Phi$. Therewith, the probability of transition $\varphi \rightarrow \psi$  
is given $Pr\{\varphi\left.\right|\psi\} =\left| (\psi,\varphi)\right|^{2}$.
Also, it was assumed that the linear operator $A$ is defined on the elements 
of linear manifold ${\cal D}(A)$ of ${\cal H}$ taking the values in
${\cal H}$. The ${\cal D}(A)$ is an everywhere closed region of definition
of $A$, namely the closure of ${\cal D}(A)$ over the norm given in
${\cal H}$ coincides with ${\cal H}$. Meanwhile, the ${\cal D}(A)$ was 
included in ${\cal D}(A^{*})$ and  $A$ coincides with the reduction of 
$A^{*}$ on ${\cal D}(A)$, because ${\cal D}(A)$ is a symmetric operator,
where the linear operator $A^{*}$ is maximal conjugated to $A$. That is,
any operator $A'$ conjugated to $A$ - $(\Psi, A'\Phi)=(A'\Psi, \Phi)$ 
at all $\Phi \in {\cal D}(A)$ and $\Psi \in {\cal D}(A')$ coincides 
with the reduction of $A^{*}$ on some linear manifold ${\cal D}(A')$
included in ${\cal D}(A^{*})$. Thus, the operator $A^{**}$ is closed 
symmetric expansion of operator $A$, namely it is a closure of $A$.
Self-conjugated operator $A$ (the closure of which self-conjugated)
allows only one self-conjugated expansion $A^{**}$. Thus, 
self-conjugated closure $\hat{N}$ of operator 
$\S_{i=1}^{\infty}
\hat{b}^{*}(f_{i})\hat{b}(f_{i})$, where $\{f_{i}\left.\right| i=1,\ldots, 
n\}$ is an arbitrary orthogonal basis in ${\cal H}^{(1)}$, was regarded
as the operator of occupation number. For the vector $\chi^{0} \in 
{}^{A}{\cal F}$ and $\chi^{0(n)}=\delta_{0n}$ one gets $<\chi^{0(n)},
\hat{N}(f)>=0$ for all $f \in {\cal H}^{(1)}$.
So, the $\chi^{0}$ is a vector of vacuum state without any particle:
$\hat{b}(f)\chi^{0}=0$ for all $f \in {\cal H}^{(1)}$. If
$f =\{ f_{i}\left.\right| i=1,2,\ldots\}$ is an arbitrary orthogonal basis 
in ${\cal H}^{(1)}$, then due to irreducibility of operators 
$\hat{b}^{*}(f_{i})\left.\right| f_{i} \in f$, the ${}^{A}{\cal H}$
includes the $0$ and whole space ${}^{A}{\cal H}$ as invariant sub-spaces 
with respect to all $\hat{b}^{*}(f)$. \\
In the sequel, to construct the 12-dimensional operator-manifold 
$\hat{G}(2.2.3)$ we consider a set $\hat{{\cal F}}$ of all sequences
$\hat{\Phi}=\{\hat{\Phi}^{(0)},\hat{\Phi}^{(1)},\ldots,\hat{\Phi}^{(n)},
\ldots \}$ with a finite number of non-zero elements, provided
\begin{equation}
\label{R329}
\begin{array}{l}
\hat{\Phi}^{(n)}_{(r_{1},\ldots,r_{n})}=\hat{\Phi}^{(1)}_{r_{1}}
\otimes \cdots \otimes \hat{\Phi}^{(1)}_{r_{n}}\in \hat{G}^{(n)},\quad
\hat{\Phi}^{(1)}_{r_{i}}=\hat{\zeta}_{r_{i}}\Phi^{(1)}_{r_{i}}\in
\hat{G}^{(1)}_{i}=\hat{\cal U}^{(1)}_{i}\otimes {\cal H}^{(1)}_{i},\\
\hat{\zeta}_{r_{i}}\equiv \S^{3}_{\alpha_{i}=1}
{\hat{\gamma}}^{r_{i}}_{(\lambda_{i},\mu_{i},\alpha_{i})}
\zeta_{r_{i}}^{(\lambda_{i},\mu_{i},\alpha_{i})}\in
\hat{\cal U}^{(1)}_{r_{i}}, \quad \hat{G}^{(n)}=\hat{\cal U}^{(n)}
\otimes \bar{\cal H}^{(n)}, \\
\hat{\cal U}^{(n)}_{(r_{1},\ldots,r_{n})}=\hat{\cal U}^{(1)}_{r_{1}},
\otimes\cdots\otimes \hat{\cal U}^{(1)}_{r_{n}}.
\end{array}
\end{equation}
Then, by analogy with eq.(3.18), the operator manifold $\hat{G}(2.2.3)$
ensued
\begin{equation}
\label{R330}
\hat{G}(2.2.3)=\S^{\infty}_{n=0}\hat{G}^{(n)}=
\S^{\infty}_{n=0}\left(\hat{\cal U}^{(n)}\otimes {\bar{\cal H}}^{(n)}\right),
\end{equation}
which decomposed as follows:
$\hat{G}(2.2.3)={}^{*}\hat{\bf R}^{22}\otimes\hat{\bf R}^{3}$,
where the linear unit operator bi-pseudo vectors 
$\{ {\hat{O}}^{r_{1}r_{2}}_{\lambda,\mu} \equiv
{\hat{O}}^{r_{1}}_{\lambda}\otimes{\hat{O}}^{r_{2}}_{\mu}\}$
served  as the basis for operator-vectors of tangent sections of  
the $2\times2$-dimensional linear 
bi-pseudo operator-space ${}^{*}\hat{\bf R}^{22}$, and
$\hat{\bf R}^{3}$ is the three-dimensional real linear operator-space with
the basis for tangent operator-vectors consisting of the ordinary 
unit operator-vectors  $\{ {\hat{\sigma}}^{r}_{\alpha}\}$. \\
At last it is useful to bring a rigorous definition of secondary
quantized form of one-particle observable $A$ in $\cal H$. Following
to [8], let consider a set of identical samples $\hat{\cal H}_{i}$ of
one-particle space ${\cal H}^{(1)}$ and operators $A_{i}$ acting in them.
To each closed linear operator $A^{(1)}$ in ${\cal H}^{(1)}$ with the
everywhere closed region of definition ${\cal D}(A^{(1)})$ following
operators correspond
\begin{equation}
\label{R331}
\begin{array}{l}
A^{(n)}_{1}=A^{(1)}\otimes I \otimes\cdots \otimes I,\\
\ldots\ldots\ldots\ldots\ldots\ldots\ldots\ldots \\
\ldots\ldots\ldots\ldots\ldots\ldots\ldots\ldots \\
A^{(n)}_{n}=I\otimes I \otimes\cdots \otimes A^{(1)}.
\end{array}
\end{equation}
Their sum $\S^{n}_{j=1}A^{(n)}_{j}$ is given on the intersection of regions of
definition of operator-terms including a linear manifold
${\cal D}(A^{(1)})\otimes\cdots \otimes {\cal D}(A^{(n)})$ being close
in $\hat{\cal H}^{(n)}$. While, the $A^{(n)}$ is a minimal closed expansion 
of this sum with ${\cal D}(A^{(n)})$. One considered linear manifold
${\cal D}(\Omega(A))$ in ${\cal H}=\S^{\infty}_{n=0}\hat{\cal H}^{(n)}$
defined as a set of all vectors $\Psi \in {\cal H}$ such as 
$\Psi^{(n)} \in {\cal D}(A^{(n)})$ and $\S^{\infty}_{n=0}\left|
A^{(n)}\Psi^{(n)}\right|^{2} < \infty$ . The manifold
${\cal D}(\Omega(A))$ is closed in ${\cal H}$. On this manifold one defines 
a closed linear operator $\Omega(A)$ acting as $\Omega(A)^{(n)}=
A^{(n)}\Psi^{(n)}$, that is $\Omega(A)\Phi = \S^{\infty}_{n=0}
A^{(n)}\Psi^{(n)}$, while the $\Omega(A)$ is self-conjugated  operator
with everywhere closed region of definition. We think of vector
$\Phi^{(n)} \in {\cal H}^{(n)}$ being in the form eq.(3.11),
where $\Phi_{i}\in {\cal D}(A)$, then $<\varphi^{(n)};A^{(n)}>=
\S^{n}_{i=1}<\varphi_{i};A>$, which allows an expansion by 
continuing onto ${\cal D}(A)$. Thus, the $A^{(n)}$ is $n$-particle 
observable corresponding to one-particle observable $A$. So. one
gets $<\varphi;\Omega(A)>=
\S^{\infty}_{n=0}<\varphi^{(n)};A^{(n)}>$ for any $\Phi_{i}\in 
{\cal D}\left(\Omega(A)\right)$. While, the $\Omega(A)$ reflects
${}^{A}{\cal D}={\cal D}\left(\Omega(A)\right)\frown {}^{A}{\cal H}$
into ${}^{A}{\cal H}$.
The reduction of $\Omega(A)$ on ${}^{A}{\cal H}$ is self-conjugated
in the region ${}^{A}{\cal D}$, because of the fact that ${}^{A}{\cal H}$
is closed sub-space of ${\cal H}$. Hence, the $\Omega(A)$ is a secondary
quantized form of one-particle observable $A$ in ${\cal H}$.\\
The main features of quantum field aspect will be completed below by 
constructing the explicit forms of wave state functions and calculating 
the matrix elements of concrete field operators.
The vacuum state reads eq.(2.16), with the normalization requirement
\begin{equation}
\label{R332}
<\chi^{0}(\nu'_{1},\nu'_{2},\nu'_{3},\nu'_{4})\mid
\chi^{0}(\nu_{1},\nu_{2},\nu_{3},\nu_{4})>=\prod_{i=1}^{4}
\delta_{\nu_{i}\nu'_{i}}.
\end{equation}
The state vectors
\begin{equation}
\label{R333}
\begin{array}{l}
\chi({\{n_{r}\}}^{N}_{1};{\{m_{r}\}}^{M}_{1};
{\{q_{r}\}}^{Q}_{1};{\{t_{r}\}}^{T}_{1};
{\{\nu_{r}\}}^{4}_{1})=
{(\hat{b}_{N}^{11})}^{n_{\scriptscriptstyle N}}\cdots
{(\hat{b}_{1}^{11})}^{n_{1}}\cdot\\
\cdot {(\hat{b}_{M}^{12})}^{m_{\scriptscriptstyle M}}\cdots
{(\hat{b}_{1}^{12})}^{m_{1}}\cdot
{(\hat{b}_{Q}^{21})}^{q_{\scriptscriptstyle Q}}\cdots
{(\hat{b}_{1}^{21})}^{q_{1}}\cdot
\cdot{(\hat{b}_{T}^{22})}^{t_{\scriptscriptstyle T}}\cdots
{(\hat{b}_{1}^{22})}^{t_{1}}
\chi^{0}(\nu_{1},\nu_{2},\nu_{3},\nu_{4}),
\end{array}
\end{equation}
where ${\{n_{r}\}}^{N}_{1}=n_{1},\ldots,n_{N}$ and so on, are the 
eigen-functions
of modified operators. They form a whole set of orthogonal vectors
\begin{equation}
\label{R334}
\begin{array}{l}
<\chi({\{n'_{r}\}}^{N}_{1};{\{m'_{r}\}}^{M}_{1};
{\{q'_{r}\}}^{Q}_{1};{\{t'_{r}\}}^{T}_{1};
{\{\nu'_{r}\}}^{4}_{1})\mid
\chi({\{n_{r}\}}^{N}_{1};{\{m_{r}\}}^{M}_{1};
{\{q_{r}\}}^{Q}_{1};{\{t_{r}\}}^{T}_{1};
{\{\nu_{r}\}}^{4}_{1})>=\\
=\displaystyle \prod_{r=1}^{N}\delta_{n_{r}n'_{r}}\cdot
\displaystyle \prod_{r=1}^{M}\delta_{m_{r}m'_{r}}\cdot
\displaystyle \prod_{r=1}^{Q}\delta_{q_{r}q'_{r}}\cdot
\displaystyle \prod_{r=1}^{T}\delta_{t_{r}t'_{r}}\cdot
\displaystyle \prod_{r=1}^{4}\delta_{\nu_{r}\nu'_{r}}.
\end{array}
\end{equation}
Considering an arbitrary superposition
\begin{equation}
\label{R335}
\chi=
\S_{\begin{array}{l}
{\scriptstyle n_{1},...,n_{N}}=0 \\
{\scriptstyle m_{1},...,m_{M}=0} \\
{\scriptstyle q_{1},...,q_{Q}=0}  \\
{\scriptstyle t_{1},...,t_{T}=0}
\end{array}}^{1}
c'({\{n_{r}\}}^{N}_{1};{\{m_{r}\}}^{M}_{1};
{\{q_{r}\}}^{Q}_{1};{\{t_{r}\}}^{T}_{1})
\chi({\{n_{r}\}}^{N}_{1};{\{m_{r}\}}^{M}_{1};
{\{q_{r}\}}^{Q}_{1};{\{t_{r}\}}^{T}_{1};
{\{\nu_{r}\}}^{4}_{1}),
\end{equation}
the coefficients $c'$ of expansion are the corresponding amplitudes
of probabilities:
\begin{equation}
\label{R336}
<\chi\mid\chi>=
\S_{\begin{array}{l}
{\scriptstyle n_{1},...,n_{N}=0} \\
{\scriptstyle m_{1},...,m_{M}=0} \\
{\scriptstyle q_{1},...,q_{Q}=0}  \\
{\scriptstyle t_{1},...,t_{T}=0}
\end{array}}^{1}
{\left| c'({\{n_{r}\}}^{N}_{1};{\{m_{r}\}}^{M}_{1};
{\{q_{r}\}}^{Q}_{1};{\{t_{r}\}}^{T}_{1}) \right|}^{2}.
\end{equation}
Taking into account eq.(3.28), the non-vanishing matrix elements of
operators $\hat{b}^{11}_{r_{k}}$ and $\hat{b}_{11}^{r_{k}}$ read
\begin{equation}
\label{R337}
\begin{array}{l}
<\chi({\{n'_{r}\}}^{N}_{1};0;0;0;1,0,0,0)\left|   \right.
\hat{b}^{11}_{r_{k}}
\chi({\{n_{r}\}}^{N}_{1};0;0;0;1,0,0,0)>=\\
\\
=<1,1\mid \hat{b}_{11}^{r'_{1}}\cdots\hat{b}_{11}^{r'_{n}}\cdot
\hat{b}^{11}_{r_{k}}\cdot
\hat{b}^{11}_{r_{n}}\cdots\hat{b}^{11}_{r_{1}}\mid 1,1>=\\
\\
= \left\{ \begin{array}{ll}
(-1)^{n'-k'}  & \mbox{if $n_{r}=n'_{r}$ for $r\neq r_{k}$ and $n_{r_{k}}=0;n'_{r_{k}}=1$}, \\
0           & \mbox{otherwise},
\end{array}  \right.   \\
\\
<\chi({\{n'_{r}\}}^{N}_{1};0;0;0;1,0,0,0)\left|   \right.
\hat{b}_{11}^{r_{k}}
\chi({\{n_{r}\}}^{N}_{1};0;0;0;1,0,0,0)>=\\
\\
=<1,1\mid \hat{b}_{11}^{r'_{1}}\cdots\hat{b}_{11}^{r'_{n}}\cdot
\hat{b}_{11}^{r_{k}}\cdot
\hat{b}^{11}_{r_{n}}\cdots\hat{b}^{11}_{r_{1}}\mid 1,1>=\\
\\
= \left\{ \begin{array}{ll}
(-1)^{n-k}  & \mbox{if $n_{r}=n'_{r}$ for $r\neq r_{k}$ and $n'_{r_{k}}=0;n_{r_{k}}=1$}, \\
0           & \mbox{otherwise},
\end{array} \right.
\end{array}
\end{equation}
where one denoted
\begin{equation}
\label{R338}
n=\S_{r=1}^{N}n_{r}, \quad n'=\S_{r=1}^{N}n'_{r},
\end{equation}
the $r_{k}$ and $r'_{k}$ are $k$-th and $k'$-th terms of regulated sets of
$\{r_{1},\ldots ,r_{n}\} \quad (r_{1}<r_{2}<\cdots <r_{n})$ and
$\{r'_{1},\ldots ,r'_{n}\} \quad (r'_{1}<r'_{2}<\cdots <r'_{n})$,
respectively.
Continuing along this line we get a whole set of explicit forms of
matrix elements of the rest of operators
$\hat{b}^{\lambda \mu}_{r_{k}}$ and $\hat{b}_{\lambda \mu}^{r_{k}}$.
Thus 
\begin{equation}
\label{R339}
\begin{array}{l}
\S_{\{\nu_{r}\}=0}^{1}<\chi^{0}\mid\hat{\Phi}(\zeta)
\mid\chi>=\S_{\{\nu_{r}\}=0}^{1}<\chi^{0}\mid
{\hat{\gamma}}_{(\lambda,\mu,\alpha)}^{r}{\Phi}^{(\lambda,\mu,\alpha)}_{r}
(\zeta)\mid\chi>= \\
=\S_{r=1}^{N}c'_{n_{r}}e_{(1,1,\alpha)}^{n_{r}}
{\Phi}^{(1,1,\alpha)}_{n_{r}}+
\S_{r=1}^{M}c'_{m_{r}}e_{(1,2,\alpha)}^{m_{r}}
{\Phi}^{(1,2,\alpha)}_{m_{r}}+ 
\S_{r=1}^{Q}c'_{q_{r}}e_{(2,1,\alpha)}^{q_{r}}
{\Phi}^{(2,1,\alpha)}_{q_{r}}+
\S_{r=1}^{T}c'_{t_{r}}e_{(2,2,\alpha)}^{t_{r}}
{\Phi}^{(2,2,\alpha)}_{t_{r}},
\end{array}
\end{equation}
provided
\begin{equation}
\label{R340}
\begin{array}{ll}
c'_{n_{r}}\equiv\delta_{1n_{r}}c'(0,\ldots, n_{r},\ldots,0;0;0;0),\quad
c'_{m_{r}}\equiv\delta_{1m_{r}}c'(0;0,\ldots m_{r},\ldots,0;0;0),\\
c'_{q_{r}}\equiv\delta_{1q_{r}}c'(0;0;0,\ldots q_{r},\ldots,0;0),\quad
c'_{t_{r}}\equiv\delta_{1t_{r}}c'(0;0;0;0,\ldots t_{r},\ldots,0).
\end{array}
\end{equation}
Hereinafter we change the notation to
\begin{equation}
\label{R341}
\begin{array}{l}
\bar{c}(r^{11})=c'_{n_{r}},\quad \bar{c}(r^{21})=c'_{q_{r}},
\quad N_{11}=N, \quad N_{21}=Q, \\
\bar{c}(r^{12})=c'_{m_{r}}, \quad \bar{c}(r^{22})=c'_{t_{r}},
\quad N_{12}=M, \quad N_{22}=T,
\end{array}
\end{equation}
and make use of
\begin{equation}
\label{R342}
\begin{array}{ll}
F_{r^{\lambda\mu}}= \S_{\alpha}
e_{(\lambda,\mu,\alpha)}^{r^{\lambda\mu}}
{\Phi}^{(\lambda,\mu,\alpha)}_{r^{\lambda\mu}},
\quad
F^{r^{\lambda\mu}}= \S_{\alpha}
e^{(\lambda,\mu,\alpha)}_{{r}^{\lambda\mu}}
{\Phi}_{(\lambda,\mu,\alpha)}^{r^{\lambda\mu}}=\bar{F}_{r^{\lambda\mu}},\\
\S_{\{\nu_{r}\}=0}^{1}<\chi^{0}\mid\hat{A}\mid\chi>
\equiv<\chi^{0}\parallel\hat{A}\parallel\chi>,
\quad
\S_{\{\nu_{r}\}=0}^{1}<\chi\mid\hat{A}\mid\chi^{0}>
\equiv<\chi\parallel\hat{A}\parallel\chi^{0}>.
\end{array}
\end{equation}
Then, the matrix elements of operator-vector and co-vector fields
take the forms
\begin{equation}
\label{R343}
\begin{array}{ll}
<\chi^{0}\parallel\hat{\Phi}(\zeta)\parallel\chi>=
\S_{\lambda\mu=1}^{2}\S_{r^{\lambda\mu}=1}^{N_{\lambda\mu}}
\bar{c}(r^{\lambda\mu})F_{r^{\lambda\mu}}(\zeta),\\
<\chi\parallel\bar{\hat{\Phi}}(\zeta)\parallel\chi^{0}>=
\S_{\lambda\mu=1}^{2}\S_{r^{\lambda\mu}=1}^{N_{\lambda\mu}}
{\bar{c}}^{*}(r^{\lambda\mu})F^{r^{\lambda\mu}}(\zeta).
\end{array}
\end{equation}
Explicitly the matrix element of $n$-operator
vector fields reads
\begin{equation}
\label{R344}
\begin{array}{l}
\FFr{1}{\sqrt{n!}}<\chi^{0}\parallel\hat{\Phi}({\zeta}_{1})
\cdots \hat{\Phi}({\zeta}_{n})\parallel\chi>=\\
=\FFr{1}{\sqrt{n!}}\left\{\S_{\lambda\mu=1}^{2} \right\}_{1}^{n}
\S_{r_{1}^{\lambda\mu},\ldots r_{n}^{\lambda\mu}=1}^{N_{\lambda\mu}}
\bar{c}(r_{1}^{\lambda\mu},\ldots , r_{n}^{\lambda\mu})
\S_{\sigma\in S_{n}}sgn(\sigma)
F_{r_{1}^{\lambda\mu}} (\zeta_{1}) \cdots
F_{r_{n}^{\lambda\mu}} (\zeta_{n})=\\
=\FFr{1}{\sqrt{n!}}\left\{\S_{\lambda\mu=1}^{2} \right\}_{1}^{n}
\S_{r_{1}^{\lambda\mu},\ldots r_{n}^{\lambda\mu}=1}^{N_{\lambda\mu}}
\bar{c}(r_{1}^{\lambda\mu},\ldots ,r_{n}^{\lambda\mu})
\left\|{ \begin{array}{ll}
F_{r_{1}^{\lambda\mu}}(\zeta_{1})\cdots F_{r_{n}^{\lambda\mu}}(\zeta_{1})\\
\cdots\cdots\cdots\cdots\cdots\cdots \\
\cdots\cdots\cdots\cdots\cdots\cdots \\
F_{r_{1}^{\lambda\mu}}(\zeta_{n})\cdots F_{r_{n}^{\lambda\mu}}(\zeta_{n})
\end{array}} \right\|.
\end{array}
\end{equation}
where $\left\{\S_{\lambda\mu}^{2} \right\}_{1}^{n}
\equiv \S_{\lambda_{1}\mu_{1}}^{2}\ldots\S_{\lambda_{n}\mu_{n}}^{2}$
and $r^{\lambda\mu}_{i}\equiv r^{\lambda_{i}\mu_{i}}$.
Here we make use of
\begin{equation}
\label{R345}
\begin{array}{l}
\bar{c}(r^{11}_{1},\ldots,r^{11}_{n})=c'(n_{1},\ldots,n_{N};0;0;0),
\quad
\bar{c}(r^{12}_{1},\ldots,r^{12}_{n})=c'(0;m_{1},\ldots,m_{M};0;0),\\
\bar{c}(r^{21}_{1},\ldots,r^{21}_{n})=c'(0;0;q_{1},\ldots,q_{Q};0),
\quad
\bar{c}(r^{22}_{1},\ldots,r^{22}_{n})=c'(0;0;0;t_{1},\ldots,t_{T}).
\end{array}
\end{equation}
Analogous expression can be obtained for co-vector fields too.
Following to this procedure, 
the anticommutation relations ensued
\begin{equation}
\label{R346}
<\chi_{-}\mid \{ {\hb_{i}}^{+}_{r},{\hb_{i}}_{+}^{r'} \}\mid\chi_{-}>=
{\delta}^{r'}_{r}, \quad
<\chi_{+}\mid \{ {\hb_{i}}^{-}_{r},{\hb_{i}}_{-}^{r'} \}\mid\chi_{+}>=
{\delta}^{r'}_{r},
\end{equation}
provided as usual
${\hgam_{i}}^{(\lambda\alpha)}_{r}=
{\he_{i}}^{(\lambda\alpha)}_{r}{\hb_{i}}^{\lambda}_{(r\alpha)},\quad
{\hgam_{i}}_{(\lambda\alpha)}^{r}=
{\he_{i}}_{(\lambda\alpha)}^{r}{\hb_{i}}_{\lambda}^{(r\alpha)}, \quad
(r\alpha)\Rightarrow r,$
where the functions $\chi_{\pm}$ have the form eq.(2.16).
The state functions
\begin{equation}
\label{R347}
\begin{array}{l}
\chi({\{n_{r}\}}^{N}_{1};{\{m_{r}\}}^{M}_{1};
{\{q_{r}\}}^{Q}_{1};{\{t_{r}\}}^{T}_{1};
\lambda;\mu)=\\
\\
={({\hb_{\eta}}_{N}^{+})}^{n_{\scriptscriptstyle N}}\cdots
{({\hb_{\eta}}_{1}^{+})}^{n_{1}}\cdot
{({\hb_{\eta}}_{M}^{-})}^{m_{\scriptscriptstyle M}}\cdots
{({\hb_{\eta}}_{1}^{-})}^{m_{1}}\cdot
{({\hb_{u}}_{Q}^{+})}^{q_{\scriptscriptstyle Q}}\cdots
{({\hb_{u}}_{1}^{+})}^{q_{1}}\cdot \\
\cdot{({\hb_{u}}_{T}^{-})}^{t_{\scriptscriptstyle T}}\cdots
{({\hb_{u}}_{1}^{-})}^{t_{1}}\cdot
\chi_{-}(\lambda)\chi_{+}(\mu),
\end{array}
\end{equation}
are normalized eq.(3.32) and formed a whole set of orthogonal
eigen-functions of the corresponding operators of occupation numbers
${\hN_{i}}^{\lambda}_{r}={\hb_{i}}_{r}^{\lambda}{\hb_{i}}^{r}_{\lambda}.$
They have the corresponding expectation values 0,1.
Finally we introduce the function $\chi$, which is an arbitrary superposition
and  has the form eq.(3.35). Incorporating with
anticommutation relations eq.(3.46)
one gets a whole set of non-vanishing matrix elements of operators
${\hb_{i}}^{\lambda}_{r}$ and ${\hb_{i}}_{\lambda}^{r}$ .
So define 
\begin{equation}
\label{R348}
\begin{array}{l}
\S_{\lambda =1}^{2}
<\chi({\{n'_{r}\}}^{N}_{1};0;0;0;\lambda;0)\left| \right.
{\hb_{\eta}}^{+}_{r_{k}}
\chi({\{n_{r}\}}^{N}_{1};0;0;0;\lambda;0)>=\\
\\
= \left\{ \begin{array}{ll}
(-1)^{n'-k'}  & \mbox{if $n_{r}=n'_{r}$ for $r\neq r_{k}$ and $n_{r_{k}}=0;n'_{r_{k}}=1,$} \\
0           & \mbox{otherwise},
\end{array}  \right. \\
\\
\S_{\lambda =1}^{2}
<\chi({\{n'_{r}\}}^{N}_{1};0;0;0;\lambda;0)\left| \right.
{\hb_{\eta}}_{+}^{r_{k}}
\chi({\{n_{r}\}}^{N}_{1};0;0;0;\lambda;0)>=\\
\\
= \left\{ \begin{array}{ll}
(-1)^{n-k}  & \mbox{if $n_{r}=n'_{r}$ for $r\neq r_{k}$ and $n'_{r_{k}}=0;n_{r_{k}}=1,$} \\
0           & \mbox{otherwise},
\end{array}  \right.
\end{array}
\end{equation}
or
\begin{equation}
\label{R349}
\begin{array}{l}
\S_{\mu =1}^{2}
<\chi(0;{\{m'_{r}\}}^{M}_{1};0;0;0;\mu)\left| \right.
{\hb_{\eta}}^{-}_{r_{k}}
\chi(0;{\{m_{r}\}}^{M}_{1};0;0;0;\mu)>=\\
\\
= \left\{ \begin{array}{ll}
(-1)^{m'-k'}  & \mbox{if $m_{r}=m'_{r}$ for $r\neq r_{k}$ and $m_{r_{k}}=0;m'_{r_{k}}=1,$} \\
0           & \mbox{otherwise},
\end{array}  \right.\\
\\
\S_{\mu =1}^{2}
<\chi(0;{\{m_{r}\}}^{N}_{1};0;0;0;\mu)\left| \right.
{\hb_{\eta}}_{-}^{r_{k}}
\chi(0;{\{m_{r}\}}^{M}_{1};0;0;0;\mu)>=\\
\\
= \left\{ \begin{array}{ll}
(-1)^{m-k}  & \mbox{if $m_{r}=m'_{r}$ for $r\neq r_{k}$ and $m'_{r_{k}}=0;m_{r_{k}}=1,$} \\
0           & \mbox{otherwise},
\end{array}  \right.
\end{array}
\end{equation}
where
\begin{equation}
\label{R350}
n=\S_{r=1}^{N}n_{r}, \quad n'=\S_{r=1}^{N}n'_{r}, \quad
m=\S_{r=1}^{M}m_{r}, \quad m'=\S_{r=1}^{M}m'_{r},
\end{equation}
and so forth.
In used notation we get
\begin{equation}
\label{R351}
<\chi^{0}\parallel\hps_{\eta}(\eta)\parallel\chi>=
\S_{\lambda=1}^{2}\S_{r^{\lambda}}^{N_{\lambda}}
\bar{c}(r^{\lambda})
{\F_{\eta}}_{r^{\lambda}}(\eta),\quad
<\chi\parallel\bar{\hps_{\eta}}(\eta)\parallel\chi^{0}>=
\S_{\lambda=1}^{2}\S_{r^{\lambda}=1}^{N_{\lambda}}
{\bar{c}}^{*}(r^{\lambda}){\bar{\F_{\eta}}}_{r^{\lambda}}(\eta).
\end{equation}
and
\begin{equation}
\label{R352}
<\chi^{0}\parallel\hps_{u}(u)\parallel\chi>=
\S_{\lambda=1}^{2}\S_{r^{\lambda}=1}^{N_{\lambda}}
\bar{c}(r^{\lambda}){\F_{u}}_{r^{\lambda}}(u), \quad
<\chi\parallel\bar{\hps_{u}}(u)\parallel\chi^{0}>=
\S_{\lambda=1}^{2}\S_{r^{\lambda}=1}^{N_{\lambda}}
{\bar{c}}^{*}(r^{\lambda}){\bar{\F_{u}}}_{r^{\lambda}}(u).
\end{equation}
provided
\begin{equation}
\label{R353}
{\F_{\eta}}_{r^{\lambda}}(\eta)= \S_{\alpha}
{\e1_{\eta}}^{(\lambda\alpha)}_{{r}^{\lambda}}
{\ps1_{\eta}}_{(\lambda\alpha)}^{{r}^{\lambda}}(\eta),\quad
{\F_{u}}_{r^{\lambda}}(u)= \S_{\alpha}
{\e1_{u}}^{(\lambda\alpha)}_{{r}^{\lambda}}
{\ps1_{u}}_{(\lambda\alpha)}^{{r}^{\lambda}}(u),\\
\end{equation}
\begin{equation}
\label{R354}
N^{1}_{1}=N, \quad N^{1}_{2}=M, \quad
N^{2}_{1}=Q, \quad N^{2}_{2}=T.
\end{equation}
Thus
\begin{equation}
\label{R355}
\begin{array}{l}
<\chi^{0}\parallel\hps_{\eta}(\eta)\bhp_{\eta}(\eta')\parallel\chi^{0}>=
\S_{\lambda\tau=1}^{2}\S_{r^{\lambda}=1}^{N^{1}_{\lambda}}
\S_{r^{\tau}=1}^{N^{1}_{\tau}}
\bar{c}(r^{\lambda})\bar{c}^{*}(r^{\tau}){\F_{\eta}}_{r^{\lambda}}(\eta),
{\bar{\F_{\eta}}}_{r^{\tau}}(\eta'), \\
<\chi^{0}\parallel\hps_{u}(u)\bhp_{u}(u')\parallel\chi^{0}>=
\S_{\lambda\tau=1}^{2}\S_{r^{\lambda}=1}^{N^{2}_{\lambda}}
\S_{r^{\tau}=1}^{N^{2}_{\tau}}
\bar{c}(r^{\lambda})\bar{c}^{*}(r^{\tau}){\F_{u}}_{r^{\lambda}}(u),
{\bar{\F_{u}}}_{r^{\tau}}(u').
\end{array}
\end{equation}

\section {Differential Geometric Aspect}
\label {diff}
The other interesting offshoot of this generalization is a geometric aspect,
the major points of which will be discussed in this section.
The set of operators $\{ \hat{\gamma}^{r}_{(\lambda,\mu,\alpha)} \}$ is
the basis for all operator-vectors of tangent section ${\hat{\bf T}}
_{\Phi_{p}}$
of principle bundle with the base $\hat{G}(2.2.3)$ at the point 
${\bf \Phi}_{p}=\left. {\bf \Phi}({\bf \zeta}(t))\right|_{t=0} \in
\hat{G}(2.2.3)$. We assumed that the $\lambda (t):R^{1}\rightarrow G(2.2.3)$
is a maximal curve passing through the point $p=\lambda (0)\in\hat{G}(2.2.3)$. 
The smooth differentiable functions 
$\{\Phi^{(\lambda,\mu,\alpha)}({\bf \zeta}(t))\}$, 
belonging to the ring of functions of $C^{\infty}$-class, are local 
coordinates in the open neighborhood of the point ${\bf \Phi}_{p} \in 
\hat{\cal U}_{\Phi}$. Passing through ${\bf \Phi}_{p}$ 
a maximal curve has the coordinates 
$\{\Phi^{(\lambda,\mu,\alpha)}({\bf \zeta}(t))\}$.
So, the smooth field of tangent operator-vector 
$\hat{\bf A}({\bf \Phi}({\bf \zeta}))$ is a class of equivalency of the 
curves ${\bf f}({\bf \Phi}({\bf \zeta}))$,
${\bf f}({\bf \Phi}({\bf \zeta}(0)))={\bf \Phi}_{p}$.
While, the operator-differential $\hat{d}\,A^{t}_{p}$ of the flux
$A^{t}_{p}:\hat{G}(2.2.3) \rightarrow \hat{G}(2.2.3)$ at the point
${\bf \Phi}_{p}$ with the velocity fields 
$\hat{\bf A}({\bf \Phi}({\bf \zeta}))$ was defined by one-parameter
group of diffeomorphismes given for the curve 
${\bf \Phi}({\bf \zeta}(t)):R^{1} \rightarrow \hat{G}(2.2.3)$, provided
${\bf \Phi}({\bf \zeta}(0))={\bf \Phi}_{p}$ and
$\widehat{\dot{{\bf \Phi}}}({\bf \zeta}(0))=\hat{\bf A}_{p}$
\begin{equation}
\label{R41}
\hat{d}\,A^{t}_{p}({\bf A})=\left.\FFr{\hat{d}}{d\,t}\right|_{t=0}
A^{t}({\bf \Phi}({\bf \zeta}(t)))=\hat{\bf A}({\bf \Phi}({\bf \zeta}))=
\hat{\gamma}^{r}_{(\lambda,\mu,\alpha)}A_{p}^{(\lambda,\mu,\alpha)},
\end{equation}
where the $\{ A_{p}^{(\lambda,\mu,\alpha)}\}$ are the components of
$\hat{\bf A}$ in the basis 
$\{ \hat{\gamma}^{r}_{(\lambda,\mu,\alpha)} \}$. According to eq.(4.1),
in holonomic coordinate basis 
$ \hat{\gamma}^{r}_{(\lambda,\mu,\alpha)} \rightarrow \left(
\hat{\partial}\left. \right/ \partial \Phi^{(\lambda,\mu,\alpha)}_{r}
({\bf \zeta}(t))\right)_{p}$ one gets
\begin{equation}
\label{R42}
A_{p}^{(\lambda,\mu,\alpha)} =\left. \FFr{\partial \Phi^{(\lambda,\mu,\alpha)}_{r}}
{\partial \zeta^{(\tau,\nu,\beta)}_{r}}\FFr{d\,
\zeta^{(\tau,\nu,\beta)}_{r}}{d\,t}\right|_{p}.
\end{equation}
The operator-tensor $\hat{\bf T}$ of 
$\widehat{(n,0)}$-type at the point ${\bf \Phi}_{p}$ is a linear
function belonging to the space
\begin{equation}
\label{R43}
\hat{\bf T}^{n}_{0}=\underbrace{\hat{\bf T}_{\Phi_{p}}\otimes\cdots \otimes 
\hat{\bf T}_{\Phi_{p}}}_{n},
\end{equation}
where $\otimes $ stands for tensor product.
It sets up a correspondence between the element $(\hat{\bf A}_{1},\ldots,
\hat{\bf A}_{n})$ of $\hat{\bf T}^{n}_{0}$ and the number $T(\hat{\bf A}_{1},
\ldots,\hat{\bf A}_{n})$, provided by linearity 
\begin{equation}
\label{R44}
\begin{array}{l}
T(\hat{\bf A}_{1},\ldots,\alpha \hat{\bf A}_{i},\ldots,\hat{\bf A}_{n})=
\alpha T(\hat{\bf A}_{1},\ldots,\hat{\bf A}_{n}), \quad \forall \alpha \in R,
\quad \hat{\bf A}_{i} \in \hat{\bf T}_{\Phi_{p}},\\
(T_{1} + T_{2})(\hat{\bf A}_{1},\ldots,\hat{\bf A}_{n})=
T_{1}(\hat{\bf A}_{1},\ldots,\hat{\bf A}_{n}) +
T_{2}(\hat{\bf A}_{1},\ldots,\hat{\bf A}_{n}).
\end{array}
\end{equation}
In the basis $\{ \hat{\gamma}^{r}_{(\lambda,\mu,\alpha)} \}$, 
for any $\hat{\bf T} \in \hat{\bf T}^{n}$, one has 
\begin{equation}
\label{R45}
\hat{\bf T}=T^{(\lambda_{1},\mu_{1},\alpha_{1})\cdots
(\lambda_{n},\mu_{n},\alpha_{n})}
_{r_{1}\cdots r_{n}}\hat{\gamma}^{r_{1}}_{(\lambda_{1},\mu_{1},\alpha_{1})}
\otimes \cdots \otimes
\hat{\gamma}^{r_{n}}_{(\lambda_{n},\mu_{n},\alpha_{n})}
\end{equation}
To render our discussion here more transparent, below
according to eq.(3.43), an explicit form of matrix element of operator-tensor
$\hat{\bf T}=\hat{\Phi}(\zeta_{1})\otimes\cdots\otimes\hat{\Phi}(\zeta_{n})$
reads
\begin{equation}
\label{R46}
\begin{array}{l}
\FFr{1}{\sqrt{n!}}<\chi^{0}\parallel\hat{\Phi}({\zeta}_{1})\otimes
\cdots \otimes\hat{\Phi}({\zeta}_{n})\parallel\chi>=\\
=\left\{\S_{\lambda\mu=1}^{2}\right\}_{1}^{n}
\S_{r_{1}^{\lambda\mu},\ldots ,r_{n}^{\lambda\mu}=1}^{N_{\lambda\mu}}
\bar{c}(r_{1}^{\lambda\mu},\ldots ,r_{n}^{\lambda\mu})
\S_{\sigma\in S_{n}}sgn(\sigma)
F_{r_{1}^{\lambda\mu}}(\zeta_{1})\otimes\cdots\otimes
F_{r_{n}^{\lambda\mu}}(\zeta_{n})= \\
=\left\{\S_{\lambda\mu=1}^{2}\right\}_{1}^{n}
\S_{r_{1}^{\lambda\mu},\ldots ,r_{n}^{\lambda\mu}=1}^{N_{\lambda\mu}}
\bar{c}(r_{1}^{\lambda\mu},\ldots ,r_{n}^{\lambda\mu})
F_{r_{1}^{\lambda\mu}}(\zeta_{1})\wedge\cdots\wedge
F_{r_{n}^{\lambda\mu}}(\zeta_{n}),
\end{array}
\end{equation}
where $\wedge$ stands for exterior product. 
If we collect together the results just established we have finally arrived
at important conclusion, that {\em by constructing matrix 
elements of operator-tensors of $\hat{G}(2.2.3)$ one produces the 
exterior products on wave manifold $\widetilde{G}(2.2.3)$}. 
There up on, {\em the matrix elements of symmetric 
operator-tensor identically equal zero}. \\
The linear operator-form of $1$ degree $\hat{\bf \omega}^{1}$ is a linear
operator-valued function on $\hat{\bf T}_{\Phi_{p}}$, namely 
$\hat{\bf \omega}^{1}(\hat{\bf A}_{p}):\hat{\bf T}_{\Phi_{p}} 
\rightarrow \hat{R}$,
where $\hat{\bf A}_{p} \in \hat{\bf T}_{\Phi_{p}}$, and the operator
$\hat{\bf \omega}^{1}(\hat{\bf A})=<\hat{\bf \omega}^{1},{\bf A}> 
\in \hat{R}$ corresponded to $\hat{\bf A}_{p}$ at the point 
${\bf \Phi}_{p}$, provided, according to eq.(3.43), with
\begin{equation}
\label{R47}
<\chi \| \hat{\bf \omega }^{1} \|
\chi^{0} >=
\S_{\lambda,\mu=1}^{2}
\S_{r^{\lambda\mu}=1}^{N_{\lambda\mu}} \hat{c}^{*}(r^{\lambda\mu})
{\bf \omega }^{1}_{r^{\lambda\mu}},
\end{equation}
where ${\bf \omega}^{1}_{r^{\lambda\mu}}=
e_{r^{\lambda\mu}}^{(\lambda,\mu,\alpha)}
\omega^{r^{\lambda\mu}}_{(\lambda,\mu,\alpha)}$, the 
$<{\bf \omega}^{1}_{r^{\lambda\mu}},{\bf A}>=
\omega ^{1}_{r^{\lambda\mu}}({\bf A})$ is a linear form on 
${\bf T}_{p}$, and 
\begin{equation}
\label{R48}
\begin{array}{l}
\hat{\bf \omega}^{1}(\lambda_{1}\hat{\bf A}_{1} + 
\lambda_{2}\hat{\bf A}_{2})=\lambda_{1}\hat{\bf \omega}^{1}(\hat{\bf A}_{1})+
\lambda_{2}\hat{\bf \omega}^{1}(\hat{\bf A}_{2}),\\
\forall \lambda_{1},\lambda_{2} \in R,\quad
\hat{\bf A}_{1}, \hat{\bf A}_{2} \in \hat{\bf T}_{\Phi_{p}}.
\end{array}
\end{equation}
The set of all linear operator-forms at the point ${\bf \Phi}_{p}$
filled operator-vector space $\hat{\bf T}_{\Phi_{p}}^{*}$ being
dual to $\hat{\bf T}_{\Phi_{p}}$. While, the
$\{ \hat{\gamma}_{r}^{(\lambda,\mu,\alpha)} \}$ served as a basis for them.
The operator $n$-form may be defined as the exterior product of operator 
1-forms
\begin{equation}
\label{R49}
\begin{array}{l}
\hat{\bf \omega}^{n}(\hat{\bf A}_{1},\ldots,\hat{\bf A}_{n})=
\left(\hat{\bf \omega}^{1}_{1}\wedge \cdots \wedge 
\hat{\bf \omega}^{1}_{n}\right)
\left(\hat{\bf A}_{1},\ldots,\hat{\bf A}_{n}\right) =\\
\\
= \left\|
\begin{array}{lll}
\hat{\bf \omega}^{1}_{1}(\hat{\bf A}_{1}) 
\cdots\cdots 
&\hat{\bf \omega}^{1}_{n}(\hat{\bf A}_{1}\\
\vdots  &\vdots\\
\vdots  &\vdots\\
\hat{\bf \omega}^{1}_{1}(\hat{\bf A}_{n}) 
\cdots\cdots 
&\hat{\bf \omega}^{1}_{n}(\hat{\bf A}_{n})
\end{array}
\right \| .
\end{array}
\end{equation}
Here, as well as for the rest of this section, to facilitate writing, 
we abbreviate the set of indices $(\lambda_{i},\mu_{i},\alpha_{i})$
by the single symbol $i$. 
If $\{\hat{\gamma}_{i}^{r_{i}}\}$ and $\{\hat{\gamma}^{i}_{r_{i}}\}$
are dual basises respectively in $\hat{\bf T}_{\Phi_{p}}$ and
$\hat{\bf T}_{\Phi_{p}}^{*}$, then the
$\{\hat{\gamma}_{1}^{r_{1}}\otimes\cdots\otimes\hat{\gamma}_{p}^{r_{p}}
\otimes\hat{\gamma}^{1}_{s_{1}}\otimes\cdots\otimes\hat{\gamma}^{q}_{s_{q}}\}$
will be the basis in operator-space
\begin{equation}
\label{R410}
\hat{\bf T}^{p}_{q}=\underbrace{\hat{\bf T}_{\Phi_{p}}\otimes\cdots \otimes 
\hat{\bf T}_{\Phi_{p}}}_{p}
\otimes\underbrace{\hat{\bf T}^{*}_{\Phi_{p}}\otimes\cdots \otimes 
\hat{\bf T}^{*}_{\Phi_{p}}}_{q}.
\end{equation}
Any operator-tensor $\hat{\bf T} \in \hat{\bf T}^{p}_{q}({\bf \Phi}_{p})$
can be written
\begin{equation}
\label{R411}
\hat{\bf T}=T^{i_{1}\cdots i_{p}}_{j_{1}\cdots j_{q}}
\left(r_{1},\ldots,r_{p},s_{1},\ldots,s_{q}\right)
\hat{\gamma}_{i_{1}}^{r_{1}}\otimes\cdots\otimes\hat{\gamma}_{i_{p}}^{r_{p}}
\otimes\hat{\gamma}^{j_{1}}_{s_{1}}\otimes\cdots\otimes
\hat{\gamma}^{j_{q}}_{s_{q}},
\end{equation}
where 
\begin{equation}
\label{R412}
T^{i_{1}\cdots i_{p}}_{j_{1}\cdots j_{q}} 
\left(r_{1},\ldots,r_{p},s_{1},\ldots,s_{q}\right)=
T\left( 
\hat{\gamma}^{i_{1}}_{r_{1}}\otimes\cdots\otimes\hat{\gamma}^{i_{p}}_{r_{p}}
\otimes\hat{\gamma}_{j_{1}}^{s_{1}}\otimes\cdots\otimes
\hat{\gamma}_{j_{q}}^{s_{q}}
\right)
\end{equation}
are the components of $\hat{\bf T}$ in dual basises
$\{\hat{\gamma}_{i}^{r_{i}}\}$ and $\{\hat{\gamma}^{i}_{r_{i}}\}$.
While, the conventional rules of tensor-algebra hold for them.
In order to utilize the operator-tensor in its concrete expression,
below according to eq.(3.43), explicitly the matrix element of
operator-tensor of $\widehat{(p,q)}$-type may be written
\begin{equation}
\label{R413}
\begin{array}{l}
\FFr{1}{\sqrt{p!q!}}<\chi^{0}\parallel{\hat{T}}^{p}_{q}\parallel\chi^{0}>=
\FFr{1}{\sqrt{p!q!}}<\chi^{0}\parallel 
T^{i_{1}\ldots i_{p}}_{j_{1}\ldots\j_{q}}
{\hat{\gamma}}_{i_{1}}^{r_{1}}\otimes\cdots\otimes
{\hat{\gamma}}_{i_{p}}^{r_{p}}
\otimes
{\hat{\gamma}}^{j_{1}}_{s_{1}}\otimes\cdots\otimes
{\hat{\gamma}}^{j_{q}}_{s_{q}}\parallel\chi^{0}>= \\
\\
=\left\{\S_{\lambda,\mu=1}^{2}\right\}_{1}^{p}
\left\{\S_{\tau,\nu=1}^{2}\right\}_{1}^{q}
\S_{r_{1}^{\lambda\mu},\ldots ,r_{p}^{\lambda\mu}=1}^{N_{\lambda\mu}}
\S_{s_{1}^{\tau\nu},\ldots ,s_{q}^{\tau\nu}=1}^{N_{\tau\nu}}
\bar{c}(r_{1}^{\lambda\mu},\ldots ,r_{p}^{\lambda\mu})
{\bar{c}}^{*}(s_{1}^{\tau\nu},\ldots ,s_{q}^{\tau\nu})\times \\
\\
\times {T(r_{1}^{\lambda\mu},\ldots ,r_{p}^{\lambda\mu};
s_{1}^{\tau\nu},\ldots ,s_{q}^{\tau\nu})}_
{\sigma(j_{1}\ldots\j_{q})}^{\sigma(i_{1}\ldots i_{p})}
(e_{i_{1}}^{r_{1}^{\lambda\mu}}\wedge\cdots\wedge
e_{i_{p}}^{r_{p}^{\lambda\mu}})\otimes
(e^{j_{1}}_{s_{1}^{\tau\nu}}\wedge\cdots\wedge
e^{j_{q}}_{s_{q}^{\tau\nu}})= \\
\\
=\left\{\S_{\lambda,\mu=1}^{2}\right\}_{1}^{p}
\left\{\S_{\tau,\nu=1}^{2}\right\}_{1}^{q}
\S_{r_{1}^{\lambda\mu},\ldots ,r_{p}^{\lambda\mu}=1}^{N_{\lambda\mu}}
\S_{s_{1}^{\tau\nu},\ldots ,s_{q}^{\tau\nu}=1}^{N_{\tau\nu}}
\bar{c}(r_{1}^{\lambda\mu},\ldots ,r_{p}^{\lambda\mu})
{\bar{c}}^{*}(s_{1}^{\tau\nu},\ldots ,s_{q}^{\tau\nu})\times \\
\\
\times {T(r_{1}^{\lambda\mu},\ldots ,r_{p}^{\lambda\mu};
s_{1}^{\tau\nu},\ldots ,s_{q}^{\tau\nu})}_
{\sigma(j_{1}\ldots\j_{q})}^{\sigma(i_{1}\ldots i_{p})}
(d\Phi_{i_{1}}^{r_{1}^{\lambda\mu}}\wedge\cdots\wedge
d\Phi_{i_{p}}^{r_{p}^{\lambda\mu}})\otimes
(d\Phi^{j_{1}}_{s_{1}^{\tau\nu}}\wedge\cdots\wedge
d\Phi^{j_{q}}_{s_{q}^{\tau\nu}}),
\end{array}
\end{equation}
The matrix elements eq.(4.6) and eq.(4.13 ) are the geometric objects 
belonging to wave manifold $\widetilde{G}(2.2.3)$ 
exposing an antisymmetric part of tensor degree. \\
For any function $f:{\bf R}^{n}\rightarrow {\bf R}^{n}$ belonging to the
ring of functions of $C^{\infty}$-class on $\hat{G}(2.2.3)$,
according to eq.(4.1), one defines an operator-differential by means of
smooth reflection $\hat{d}\,f:\hat{\bf T}\left(\hat{G}(2.2.3)\right)
\rightarrow \hat{R}\quad \left(\hat{\bf T}\left(\hat{G}(2.2.3)\right)=
\displaystyle \bigcup_{\Phi_{p}}\hat{\bf T}_{\Phi_{p}}\right)$ as follows:
\begin{equation}
\label{R414}
<\hat{d}\,f,\hat{\bf A}>=(Af)^{\wedge},
\end{equation}
where it is denoted
\begin{equation}
\label{R415}
<\chi\| B^{\wedge}\|\chi^{0}>=\S_{\lambda,\mu=1}^{2}
\S_{r^{\lambda\mu}=1}^{N_{\lambda\mu}} \hat{c}^{*}(r^{\lambda\mu})
B(r^{\lambda\mu}).
\end{equation}
This rule for operator-differential is not without an important reason.
The argument for this conclusion is compulsory suggested by the properties
of the group of diffeomorphismes eq.(4.1). The ordinary laws regarding these 
changes apply that we ought to make use of
\begin{equation}
\label{R416}
<\chi\|\hat{d}\,f,\hat{\bf A} \|\chi^{0}>=\S_{\lambda,\mu=1}^{2}
\S_{r^{\lambda\mu}=1}^{N_{\lambda\mu}} \hat{c}^{*}(r^{\lambda\mu})
<d\,f,{\bf A}>_{r^{\lambda\mu}}=
\S_{\lambda,\mu=1}^{2}
\S_{r^{\lambda\mu}=1}^{N_{\lambda\mu}} \hat{c}^{*}(r^{\lambda\mu})
({\bf A}\,f)_{r^{\lambda\mu}},
\end{equation}
then in coordinate basis
\begin{equation}
\label{R417}
<d\, \Phi^{\widehat{\imath}}, \hat{\partial}\left.\right/ \partial 
\Phi^{j}>=\FFr{\partial \Phi^{\widehat{\imath}}}{\partial \Phi^{j}}=
\delta^{{\widehat{\imath}}}_{j},
\end{equation}
provided $d\, \Phi^{\widehat{\imath}}\equiv \hat{d}\Phi^{i}$, the
$\delta^{{\widehat{\imath}}}_{j}$ is being fashioned after the conventional 
symbol $\delta^{i}_{j}$ and best visualized as
\begin{equation}
\label{R418}
<\chi\|\delta^{{\widehat{\imath}}}_{j} \|\chi^{0}>=\S_{\lambda,\mu=1}^{2}
\S_{r^{\lambda\mu}=1}^{N_{\lambda\mu}} \hat{c}^{*}(r^{\lambda\mu})
\delta^{i}_{j}, \quad (i \leftrightarrow 
(\lambda_{i},\mu_{i},\alpha_{i})).
\end{equation}
Continuing along this line we define the differential operator
$n$-form $\left.{\bf \hat{\omega}}^{n}\right|_{\Phi_{p}}$ at the point
${\bf \Phi}_{p} \in \hat{G}(2.2.3)$ as the exterior operator $n$-form
on tangent operator space $\hat{\bf T}_{\Phi_{p}}$ of tangent
operator-vectors $\hat{\bf A}_{1},\ldots,\hat{\bf A}_{n}$.
That is, if the $\wedge \hat{\bf T}^{*}_{\Phi_{p}}
\left(\hat{G}(2.2.3)\right)$ means the exterior algebra on
$\hat{\bf T}^{*}_{\Phi_{p}}\left(\hat{G}(2.2.3)\right)$, then
operator $n$-form $\left.{\bf \hat{\omega}}^{n}\right|_{\Phi_{p}}$
is an element of $n$-th degree out of $\wedge \hat{\bf T}^{*}_{\Phi_{p}}$
depending on the point ${\bf \Phi}_{p} \in \hat{G}(2.2.3)$.
Hence ${\bf \hat{\omega}}^{n}=\displaystyle \bigcup_{\Phi_{p}}
\left.{\bf \hat{\omega}}^{n}\right|_{\Phi_{p}}$. Any differential operator
$n$-form belonging to dual operator-space
$\underbrace{\hat{\bf T}^{*}_{\Phi_{p}}\otimes\cdots \otimes 
\hat{\bf T}^{*}_{\Phi_{p}}}_{n}$ may be written in form
\begin{equation}
\label{R419}
{\bf \hat{\omega}}^{n}=\S_{i_{1}<\cdots<i_{n}}\alpha
_{i_{1} \cdots i_{n}} (\Phi)d\,\Phi^{\widehat{\imath}_{1}}\wedge\cdots\wedge
d\,\Phi^{\widehat{\imath}_{n}},
\end{equation}
provided by the smooth differentiable functions $\alpha
_{i_{1} \cdots i_{n}}(\Phi)\in C^{\infty}$ and basis
\begin{equation}
\label{R420}
d\,\Phi^{\widehat{\imath}_{1}}\wedge\cdots\wedge
d\,\Phi^{\widehat{\imath}_{n}}=
\S_{\sigma\in S_{n}}sgn(\sigma)
\gamma^{\sigma(\widehat{\imath}_{1}}\otimes\cdots\otimes
\gamma^{\widehat{\imath}_{n})}.
\end{equation}
So, any antisymmetric operator-tensor of $\widehat{(0,n)}$-type
reads 
\begin{equation}
\label{R421}
\hat{\bf T}^{*}=T_{i_{1} \cdots i_{n}}\gamma^{\widehat{\imath}_{1}}
\otimes\cdots\otimes \gamma^{\widehat{\imath}_{n}}
=\S_{i_{1}<\cdots<i_{n}}
T_{i_{1} \cdots i_{n}} d\,\Phi^{\widehat{\imath_{1}}}\wedge\cdots\wedge
d\,\Phi^{\widehat{\imath_{n}}}.
\end{equation}
Let the $\hat{\cal D}_{1}$ and $\hat{\cal D}_{2}$ are two compact 
convex parallelepipeds in oriented $n$-dimensional operator-space
$\hat{\bf R}^{n}$ and the $f:\hat{\cal D}_{1}\rightarrow \hat{\cal D}_{2}$
is differentiable reflection of interior of $\hat{\cal D}_{1}$ into
$\hat{\cal D}_{2}$ retaining an orientation, namely for any function
$\varphi \in C^{\infty}$ defined on $\hat{\cal D}_{2}$ if holds
$\varphi\circ f\in C^{\infty}$ and $f^{*}\varphi\left({\bf \Phi}_{p}\right)=
\varphi\left( f \left({\bf \Phi}_{p} \right)\right)$, where $f^{*}$
is an image of function $\varphi\left( f \left({\bf \Phi}_{p} \right)\right)$
on $\hat{\cal D}_{1}$ at the point ${\bf \Phi}_{p}$. Hence, the function
$f$ induced a linear reflection $\hat{d}\, f:\hat{\bf T}\left(
\hat{\cal D}_{1}\right)\rightarrow \hat{\bf T}\left(
\hat{\cal D}_{2}\right)$ as an operator-differential of $f$ implying
$\hat{d}\,f\left(\hat{\bf A}_{p} \right)\varphi=\hat{\bf A}_{p}
(\varphi\circ f)$ for any operator-vector $\hat{\bf A}_{p} \in
\hat{\bf T}_{\Phi_{p}}$ and for any function $\varphi \in C^{\infty}$ defined
in the neighborhood of ${\bf \Phi'}_{p}=f\left({\bf \Phi}_{p}\right)$.
If the function $f$ is given in form ${\Phi'}^{i}={\Phi'}^{i}\left(
\Phi^{(1,1,1)}_{p},\ldots,\Phi^{(2,2,3)}_{p}\right)$, where the
coordinate suffixes were only put forth in illustration of a point at
issue, and
$\hat{\bf A}_{p}=\left(A^{i}\left.\hat{\partial}\right/ \partial\Phi^{i}
\right)_{p}$, then in terms of local coordinates one gets
\begin{equation}
\label{R422}
\left(\hat{d}\,f\right)\hat{\bf A}_{p}=A^{i}\left(\FFr{\partial{\Phi'}^{j}}
{\partial\Phi^{i}}\right)_{p}
\left(\FFr{\hat{\partial}}
{\partial{\Phi'}^{j}}\right)_{p'}.
\end{equation}
So, if $f_{1}:\hat{\cal D}_{1}\rightarrow \hat{\cal D}_{2}$ and
$f_{2}:\hat{\cal D}_{2}\rightarrow \hat{\cal D}_{3}$ then
$\hat{d}\,\left(f_{2}\circ f_{1}\right)=\hat{d}\,f_{2}\circ \hat{d}\,
f_{1}$. The differentiable reflection 
$f:\hat{\cal D}_{1}\rightarrow \hat{\cal D}_{2}$ induced the reflection
$\hat{f}^{*}:\hat{\bf T}^{*}\left(\hat{\cal D}_{2}\right)\rightarrow 
\hat{\bf T}^{*}\left(\hat{\cal D}_{1}\right)$ conjugated to $\hat{f}_{*}$.
The latter is an operator differential of $f$, while
\begin{equation}
\label{R423}
<\hat{f}^{*}\hat{\omega'}^{1},\hat{\bf A}>_{\Phi_{p}}=
\left.<\hat{\omega'}^{1},\hat{f}_{*}\hat{\bf A}>\right|_{f\left(
\Phi_{p}\right)},
\end{equation}
where $\left.\hat{\bf A}\right|_{f\left(\Phi_{p}\right)}=
\left(\hat{d}\,f\right)\hat{\bf A}_{p}$ and $\hat{\omega'}^{1}\in 
\left.\hat{\bf T}^{*}\right|_{f\left(\Phi_{p}\right)}$. Hence
\begin{equation}
\label{R424}
\hat{f}^{*}\left(\hat{d}\,\varphi\right)=\hat{d}\,
\left(\hat{f}^{*}\varphi\right)
\end{equation}
and
\begin{equation}
\label{R425}
\begin{array}{l}
\hat{f}_{*}:\hat{\bf T}\in \hat{\bf T}^{0}_{n}\left(\Phi_{p}\right)
\rightarrow \hat{f}_{*}\hat{\bf T}\in 
\hat{\bf T}_{0}^{n}\left(f\left(\Phi_{p}\right)\right),\\
\\
\hat{f}^{*}:\hat{\bf T}\in \hat{\bf T}_{0}^{n}
\left(f\left(\Phi_{p}\right)\right)
\rightarrow \hat{f}^{*}\hat{\bf T}\in 
\hat{\bf T}^{0}_{n}\left(\Phi_{p}\right).
\end{array}
\end{equation}
That is
\begin{equation}
\label{R426}
\begin{array}{l}
\hat{f}^{*}\left.T\left(\hat{\bf A}_{1},\ldots,
\hat{\bf A}_{n}\right)\right|_{\Phi_{p}}=
\left.T\left(\hat{f}_{*}\hat{\bf A}_{1},\ldots,\hat{f}_{*}
\hat{\bf A}_{n}\right)\right|_{f\left(\Phi_{p}\right)},\\
\\
\left.T\left(\hat{f}^{*}\hat{\bf \omega}^{1}_{1},\ldots,
\hat{f}^{*}\hat{\bf \omega}^{1}_{n}\right)\right|_{\Phi_{p}}=
\hat{f}_{*}\left.T\left(\hat{\bf \omega}^{1}_{1},\ldots,\hat{f}^{*}
\hat{\bf \omega}^{1}_{n}\right)\right|
_{f\left(\Phi_{p}\right)}.
\end{array}
\end{equation}
So, for any differential operator $n$-form $\hat{\bf \omega}^{n}$ on
$\hat{\cal D}_{2}$ the reflection $f$ induced the operator $n$-form
$\hat{f}^{*}\hat{\bf \omega}^{n}$ on $\hat{\cal D}_{1}$
\begin{equation}
\label{R427}
\left(\hat{f}^{*}\hat{\bf \omega}^{n}\right)
\left.\left(\hat{\bf A}_{1},\ldots,
\hat{\bf A}_{n}\right)\right|_{\Phi_{p}}=\hat{f}_{*}\hat{\bf \omega}^{n}
\left.\left(\hat{f}_{*}\hat{\bf A}_{1},\ldots,\hat{f}_{*}
\hat{\bf A}_{n}\right)\right|_{f\left(\Phi_{p}\right)}.
\end{equation}
If $\hat{\bf \omega'}^{1}=\alpha'_{i}d\,{\Phi'}^{\widehat{\imath}}$
then
\begin{equation}
\label{R428}
\hat{f}^{*}\left(\alpha'_{i}d\,{\Phi'}^{\widehat{\imath}}\right)=
\alpha'_{i}\FFr{\partial {\Phi'}^{i}}{\partial {\Phi}^{j}}
d\,\Phi^{\widehat{\jmath}}.
\end{equation}
This can be extended to
$\hat{\bf \omega'}^{n} \rightarrow \hat{\bf \omega}^{n}$
\begin{equation}
\label{R429}
\begin{array}{l}
\hat{f}^{*}\left(
\S_{i_{1}<\cdots<i_{n}}
{T'}_{i_{1} \cdots i_{n}} d\,{\Phi'}^{\widehat{\imath_{1}}}\wedge\cdots\wedge
d\,{\Phi'}^{\widehat{\imath_{n}}}\right)=\\
=\S_{\begin{array}{l}
{\scriptstyle i_{1}<\cdots<i_{n}} \\
{\scriptstyle j_{1}<\cdots<j_{n}}
\end{array}}
{T'}_{i_{1} \cdots i_{n}} 
\FFr{\partial {\Phi'}^{i_{1}}}{\partial {\Phi}^{j^{1}}}\cdots
\FFr{\partial {\Phi'}^{i_{n}}}{\partial {\Phi}^{j^{n}}}
d\,{\Phi}^{\widehat{\imath_{1}}}\wedge\cdots\wedge
d\,{\Phi}^{\widehat{\imath_{n}}},
\end{array}
\end{equation}
namely
\begin{equation}
\label{R430}
\hat{f}^{*}\hat{\bf \omega'}^{n} = J_{\Phi}\hat{\bf \omega}^{n}=
\left( det\,d f\right)\hat{\bf \omega}^{n},
\end{equation}
where $J_{\Phi}$ is the Jacobian of reflection $J_{\Phi}=\left\|
\FFr{\partial {\Phi'}^{i}}{\partial {\Phi}^{j}}\right\|$.
While, the following relations hold
\begin{equation}
\label{R431}
\left(\hat{f}_{1}\circ\hat{f}_{2}\right)^{*}=
\hat{f}^{*}_{1}\circ\hat{f}^{*}_{2},\quad
\hat{f}^{*}\left(\hat{\bf \omega}_{1}\wedge
\hat{\bf \omega}_{2}\right)=
\hat{f}^{*}\left(\hat{\bf \omega}_{1}\right)\wedge
\hat{f}^{*}\left(\hat{\bf \omega}_{2}\right).
\end{equation}
We may consider the integration of operator $n$-form implying
\begin{equation}
\label{R432}
\IIn_{\hat{\cal D}_{1}}\hat{f}^{*}\hat{\bf \omega}^{n} = 
\IIn_{\hat{\cal D}_{2}}\hat{\bf \omega}^{n}.
\end{equation}
In general, let the $\hat{\cal D}_{1}$ is a limited convex
$n$-dimensional parallelepiped in $n$-dimensional operator space
$\hat{\bf R}^{n}$. One defined the $n$-dimensional $i$-th piece of integration
path $\hat{\sigma}^{i}$ in $\hat{G}(2.2.3)$ as $\hat{\sigma}^{i}=
\left(\hat{\cal D}_{i}, f_{i}, Or_{i}\right)$, where $\hat{\cal D}_{i}
\in \hat{\bf R}^{n}, \quad f_{i}:\hat{\cal D}_{i}\rightarrow \hat{G}(2.2.3)$
and the $Or_{i}$ is an orientation of $\hat{\bf R}^{n}$. Then, the integral
over the operator $n$-form $\hat{\bf \omega}^{n}$ along the operator
$n$-dimensional chain $\hat{c}_{n}=\S m_{i}\hat{\sigma}^{i}$ may be
written
\begin{equation}
\label{R433}
\IIn_{\hat{c}_{n}}\hat{\bf \omega}^{n}=\S m_{i}
\IIn_{\hat{\sigma}^{i}}\hat{\bf \omega}^{n}=\S m_{i}
\IIn_{\hat{\cal D}_{i}}\hat{f}^{*}\hat{\bf \omega}^{n},
\end{equation}
where the $m_{i}$ is a multiple number. Taking into account the eq.(4.7),
its matrix element yields
\begin{equation}
\label{R434}
<\chi \| \IIn_{\hat{c}_{n}}\hat{\bf \omega}^{n}\|
\chi^{0} >\rightarrow
\left\{\S_{\lambda\mu=1}^{2}\right\}_{1}^{n}
\S_{r_{1}^{\lambda\mu},\ldots ,r_{n}^{\lambda\mu}=1}^{N_{\lambda\mu}}
\S m_{i}\bar{c}(r_{1}^{\lambda\mu},\ldots ,r_{n}^{\lambda\mu})
\IIn_{\hat{\cal D}_{i}}\hat{f}^{*}\hat{\bf \omega}^{n}
(r_{1}^{\lambda\mu},\ldots ,r_{n}^{\lambda\mu}).
\end{equation}
Next we employ the analog of exterior differentiation adjusted to fit 
a discussing formalism. So, we may define the value of the 
operator $(n+1)$-form $\hat{d}\,\hat{\bf \omega}^{n}$ on $(n+1)$
operator-vectors $\hat{\bf A}_{1},\ldots,
\hat{\bf A}_{n+1}\in \hat{\bf T}_{\Phi_{p}}$ by considering 
diffeomorphic reflection $f$ of the neighborhood of the point
$0$ in $\hat{\bf R}^{n}$ into neighborhood of the point ${\bf \Phi}_{p}$
in $\hat{G}(2.2.3)$. The prototypes of operator-vectors
$\hat{\bf A}_{1},\ldots,
\hat{\bf A}_{n+1}\in \hat{\bf T}_{\Phi_{p}}\left(\hat{G}(2.2.3)\right)$
at the operator differential of $f$ belong to tangent operator space
$\hat{\bf R}^{n}$ in $0$. Then, the prototypes are the operator-vectors
$\hat{\bf \xi}_{1},\ldots,\hat{\bf \xi}_{n+1} \in \hat{\bf R}^{n}$. Let
$f$ reflects the parallelepiped $\hat{\bf \Pi}^{*}$, stretched over 
the $\hat{\bf \xi}_{1},\ldots,\hat{\bf \xi}_{n+1}$, into $(n+1)$-
dimensional piece $\hat{\bf \Pi}$ on $\hat{G}(2.2.3)$.
While the border of $n$-dimensional chain $\partial \hat{\bf \Pi}$ 
in $\hat{\bf R}^{n+1}$ defined as follows: the pieces $\hat{\sigma}^{i}$
of the chain $\partial \hat{\bf \Pi}$ are $n$-dimensional facets
$\partial \hat{\bf \Pi}_{i}$ of parallelepiped $\partial \hat{\bf \Pi}$
with the reflections $f_{i}:\hat{\bf \Pi}_{i} \rightarrow
\hat{\bf R}^{n+1}$ of embedding of the facets into $\hat{\bf R}^{n+1}$, and the
orientations $Or_{i}$ defined $\partial \hat{\bf \Pi}=\S\hat{\sigma}^{i},
\quad \hat{\sigma}^{i}=\left(\hat{\bf \Pi}_{i},f_{i},Or_{i}\right)$
Considering the curvilinear parallelepiped
\begin{equation}
\label{R435}
F\left(\hat{\bf A}_{1},\ldots,\hat{\bf A}_{n}\right)=
\IIn_{\partial \hat{\bf \Pi}}\hat{\bf \omega}^{n},
\end{equation}
one may state that the unique operator of $(n+1)$-form $\hat{\Omega}$
exists on $\hat{\bf T}_{\Phi_{p}}$, which is the major $(n+1)$-linear
part in $0$ of integral over the border of 
$F\left(\hat{\bf A}_{1},\ldots,\hat{\bf A}_{n}\right)$, namely
\begin{equation}
\label{R436}
F\left(\varepsilon\hat{\bf A}_{1},\ldots,
\varepsilon\hat{\bf A}_{n}\right)=
\varepsilon^{n+1}\hat{\Omega}
\left(\hat{\bf A}_{1},\ldots,\hat{\bf A}_{n+1}\right)+
O\left(\varepsilon^{n+1}\right),
\end{equation}
where the $\hat{\Omega}$ is independent of the choice of coordinates
used in definition of $F$. 
A prove of it may be readily furnished in analogy with the same statement
of corresponding theorem of differential geometry [9].
If in local coordinates 
\begin{equation}
\label{R437}
\hat{\bf \omega}^{n}=\S_{i_{1}<\cdots<i_{n}}
T_{i_{1} \cdots i_{n}} d\,\Phi^{\widehat{\imath_{1}}}\wedge\cdots\wedge
d\,\Phi^{\widehat{\imath_{n}}},
\end{equation}
then
\begin{equation}
\label{R438}
\hat{\Omega}=\hat{d}\,\hat{\bf \omega}^{n}=
\S_{i_{1}<\cdots<i_{n}}
\hat{d}\,T_{i_{1} \cdots i_{n}} d\,\Phi^{\widehat{\imath_{1}}}
\wedge\cdots\wedge d\,\Phi^{\widehat{\imath_{n}}}.
\end{equation}
The operator of exterior differential $\hat{d}$ commutes with the
reflection $f:\hat{G}(2.2.3)\rightarrow \hat{G}(2.2.3)$
\begin{equation}
\label{R439}
\hat{d}\,\left(\hat{f}^{*}\hat{\bf \omega}^{n}\right)=
\hat{f}^{*}\left(\hat{d}\,\hat{\bf \omega}^{n}\right).
\end{equation}
So define the exterior differential by operator-form $d\hat{\omega}$
of $(n+1)$ degree
\begin{equation}
\label{R440}
\begin{array}{l}
\hat{d}\hat{\omega}^{n}=\S_{\begin{array}{l}
{\scriptstyle i_{0}} \\
{\scriptstyle i_{1}<\ldots<i_{n}}
\end{array}}
\FFr{\partial T_{i_{1}\ldots i_{n}}}{\partial\Phi^{i_{0}}}
d\Phi^{\widehat{i_{0}}}\wedge
d\Phi^{\widehat{i_{1}}}\wedge\cdots\wedge d\Phi^{\widehat{i_{n}}}= \\
\\
=\S_{i_{1}<\ldots<i_{n}}(\hat{d}T_{i_{1}\ldots i_{n}})\wedge
d\Phi^{\widehat{i_{1}}}\wedge\cdots\wedge d\Phi^{\widehat{i_{n}}},
\end{array}
\end{equation}
then
\begin{equation}
\label{R441}
\begin{array}{l}
<\chi\parallel \hat{d}\hat{\omega}^{n} \parallel\chi^{0}>\rightarrow \\
\\
\S_{ \begin{array}{l}
{\scriptstyle i_{0}} \\
{\scriptstyle i_{1}<\ldots<i_{n}}
\end{array}}
\left\{\S_{\lambda,\mu=1}^{2}\right\}_{0}^{n}
\S_{r_{1}^{\lambda\mu},\ldots ,r_{n}^{\lambda\mu}=1}^{N_{\lambda\mu}}
\bar{c}(r_{0}^{\lambda\mu},r_{1}^{\lambda\mu},\ldots ,r_{n}^{\lambda\mu})
\times\\
\times \FFr{\partial T
(r_{1}^{\lambda\mu},\ldots ,
r_{n}^{\lambda\mu})_{i_{1}\ldots i_{n}}}
{\partial\Phi^{i_{0}}_{r_{0}^{\lambda\mu}}}
d\Phi^{i_{0}}_{r_{0}^{\lambda\mu}}\wedge
d\Phi^{i_{1}}_{r_{1}^{\lambda\mu}}\wedge\cdots\wedge 
d\Phi^{i_{n}}_{r_{n}^{\lambda\mu}}= \\
\\
=\S_{i_{1}<\ldots<i_{n}}
\left\{\S_{\lambda,\mu=1}^{2}\right\}_{1}^{n}
\S_{r_{1}^{\lambda\mu},\ldots ,r_{n}^{\lambda\mu}=1}^{N_{\lambda\mu}}
\bar{c}(r_{1}^{\lambda\mu},\ldots r_{n}^{\lambda\mu})
{(dT(r_{1}^{\lambda\mu},\ldots ,r_{n}^{\lambda\mu})}_{i_{1}\ldots i_{n}})
\wedge 
d\Phi^{i_{1}}_{r_{1}^{\lambda\mu}}\wedge\cdots\wedge 
{d\Phi}^{i_{n}}_{r_{n}^{\lambda\mu}}.
\end{array}
\end{equation}
In generalized sense, we may draw a conclusion that
{\em the matrix element 
of any geometric object of operator 
manifold $\hat{G}(2.2.3)$ yields corresponding geometric object of wave 
manifold $\widetilde{G}(2.2.3)$}. Thus, {\em all geometric objects 
belonging to the latter can be constructed by means of matrix elements of 
corresponding geometric objects of the operator manifold $\hat{G}(2.2.3)$}.\\
One final observation is worth recording. We may
consider a linear differential operator $\hat{L}_{\hat{\bf A}}$,
namely the differentiation along the direction of operator-vector 
$\hat{\bf A}$. For any function $\varphi \in C^{\infty}, \quad
\varphi:\hat{G}(2.2.3) \rightarrow R$ a derivative     along a direction of
$\hat{\bf A}$ is an operator-function $\hat{L}_{\hat{\bf A}}\varphi$
taking at the point ${\bf \Phi}_{p}$ the value
\begin{equation}
\label{R442}
\left(\hat{L}_{\hat{\bf A}}\varphi({\bf \Phi}) \right)=
\left.\FFr{\hat{d}}{d\,t}\right|_{t=0}\varphi\left(A^{t}{\bf \Phi}\right).
\end{equation}
In local coordinates $\{ \Phi^{i}\}$ the flux $A^{t}$ is given by equation
$\widehat{\dot{\Phi}}^{i}(\zeta)=A^{\widehat{\imath}}(\zeta)$. So
\begin{equation}
\label{R443}
\hat{L}_{\hat{\bf A}}=
A^{i}\FFr{\hat{\partial}}{\partial \Phi^{i}}.
\end{equation}
Let it be ${\hat{\bf A}},{\hat{\bf B}}\in
\hat{\bf T}_{\Phi_{p}}\left(\hat{G}(2.2.3)\right)$.
Then, two fluxes $A^{t}$ and $B^{t}$ are commuting only if Poison bracket
$[{\hat{\bf A}},{\hat{\bf B}}]$ is zero, while for the field
${\hat{\bf C}}=[{\hat{\bf A}},{\hat{\bf B}}]$  it holds
$\hat{L}_{\hat{\bf C}}=\hat{L}_{\hat{\bf B}}\hat{L}_{\hat{\bf A}}-
\hat{L}_{\hat{\bf A}}\hat{L}_{\hat{\bf B}}$. We may extend this up to
differentiation of operator-tensor $\hat{\bf T}$ at the point
${\bf \Phi}_{p}$ along the direction of $\hat{\bf A}$
\begin{equation}
\label{R444}
\left.\hat{L}_{\hat{\bf A}}\hat{\bf T}\right|_{\Phi_{p}}=
\displaystyle \lim_{t\rightarrow 0}\FFr{1}{t}\left\{ \left.\hat{\bf T}
\right|_{\Phi_{p}} -\hat{f}_{t^{*}} \left.\hat{\bf T}
\right|_{\Phi_{p}}\right\}.
\end{equation}
We hope that we have made good
headway by presenting a reasonable analysis of mathematical structure of
the method of operator manifold $\hat{G}(2.2.3)$.

\section {Discussion and Conclusions}
At this point we cut short our exposition of the
theory and reflect upon the results far obtained.
A number of conclusions may be drawn and the main features of suggested
theory  are outlined below. The method of operator manifold is a still wider 
generalization of secondary quantization with appropriate expansion
over the geometric objects, leading to the quantization of geometry differed 
in principle from all earlier suggested schemes. It has two important aspects
of quantum field theory and differential geometry.
We first deal with a substitution of the basis vectors of tangent section of
the manifold $G(2.2.3)$ by corresponding operators of creation and 
annihilation of quanta of geometry acting in the configuration space of
occupation numbers. While, as it was shown, the quantum of geometry may 
be regarded as a fermion field described by the theory being in close 
analogy to Dirac's conventional wave-mechanical theory of fermions with a 
spin $\vec{\FFr{1}{2}}$ treated in terms of manifold $G(2.2.3)$.
The final formulation of the method is mainly based on configuration space 
mechanics with antisymmetric state functions incorporated with geometric 
properties of corresponding elements. In view of this, one has reached
to well-grounded rigorous definition of concept of operator manifold
$\hat{G}(2.2.3)$. Completing the quantum field aspect of operator manifold
$\hat{G}(2.2.3)$, we have constructed the explicit forms of wave state 
functions and calculated the matrix elements of concrete field operators.
The differential geometric aspect, which is a
second interesting offshoot of the method of operator manifold, is
discussed in detail in last section. We defined the operator-tensors
and have drawn an important conclusion, that the matrix elements of 
operator-tensors of $\hat{G}(2.2.3)$ yield the external product on wave 
manifold $\widetilde{G}(2.2.3)$. One considered also the differential operator 
forms and employed their integration as well as exterior differentiation,
the differentiation along the direction of operator-vector, adjusted to
fit the formalism of operator manifold. One may state that
the matrix element of any geometric object of operator 
manifold $\hat{G}(2.2.3)$ yields corresponding geometric object of wave 
manifold $\widetilde{G}(2.2.3)$.

\centerline {\bf\large Acknowledgements}
\vskip 0.1\baselineskip
\noindent
I am greatly indebted to A.M.Vardanian and K.L.Yerknapetian
for support.


\begin {thebibliography}{99}
\bibitem{A3} Ter-Kazarian G.T., Astrophys. and Space Sci., 1996, v.241, 
No2, 161, hep-th/9710089.
\bibitem{\Ter86} Ter-Kazarian G.T., Selected Questions of Theoretical 
and Mathematical Physics, 1986, VINITI, N5322-B86, Moscow.
\bibitem{\Tercom} Ter-Kazarian G.T.,Comm. Byurakan Obs., 1989,v.62, 1.
\bibitem{\Ter92} Ter-Kazarian G.T., Astrophys. and Space Sci., 1992, v.194, 1.
\bibitem{A1} Ter-Kazarian G.T.,IC/94/290, ICTP Preprint, Trieste, 
Italy, 1, 1994.
\bibitem{A2} Ter-Kazarian G.T., hep-th/9510110, 1995.
\bibitem{A13} Bjorken J.D., Drell S.D., Relativistic Quantum Fields,
\indent Mc Graw-Hill, New York, 1964.
\bibitem{A8} Cook J.M., Trans. Amer. Math. Soc., 1953, 74, 222.
\bibitem{A9} Arnold V.I., 1979, Mathemathical Methods of Classical
Mechanics, Nauka, Moscow.
\end {thebibliography}
\end{document}